\documentclass{conm-p-l}

\newtheorem{theorem}{Theorem}[section]

\theoremstyle{definition}

\theoremstyle{remark}
\newtheorem{remarkplain}[theorem]{Remark}
\newtheorem{exampleplain}[theorem]{Example}

\numberwithin{equation}{section}

\usepackage{amssymb}
\usepackage{amsmath}
\usepackage{tabularx}
\usepackage[dvips]{graphicx}
\usepackage[arrow,matrix,tips,frame,color,line,poly]{xy}

\newcommand{\Aut}{\text{Aut}}
\newcommand{\diag}{\text{diag}}

\newcommand{\Out}{\text{Out}}

\renewcommand{\int}{\text{int}}

\newcommand{\End}{\text{End}}

\newcommand{\Norm}{\text{Norm}}

\newcommand{\X}{\mathcal X}
\newcommand{\Xr}{\mathcal X^r}

\newcommand{\caZ}{\mathcal Z}
\newcommand{\caZr}{\mathcal Z^r}

\newcommand{\I}{\mathcal I}

\newcommand{\tX}{\widetilde{\mathcal X}}
\newcommand{\tXr}{\widetilde{\mathcal X}^r}

\newcommand{\R}{\mathbb R}
\newcommand{\C}{\mathbb C}
\newcommand{\Z}{\mathbb Z}

\newcommand{\Ztoo}{\mathbb Z_2}

\newcommand{\Q}{\mathbb Q}

\newcommand{\caB}{\mathcal B}

\newcommand{\G}{G}

\newcommand{\h}{\mathfrak h}

\newcommand{\ch}[1]{#1^\vee}
\newcommand{\chG}{\ch{G}}

\renewcommand{\sec}[1]{\section{#1}
\renewcommand{\theequation}{\thesection.\arabic{equation}}
  \setcounter{equation}{0}}
\newcommand{\subsec}[1]{\subsection{#1}
\renewcommand{\theequation}{\thesection.\arabic{equation}}}

\newcommand{\g}{\mathfrak g}
\newcommand{\GGamma}{G^\Gamma}

\newcommand\inv{^{-1}}
\newcommand\bs{\backslash}

\begin{document}

\author{Jeffrey Adams}

\address{Department of Mathematics,
University of Maryland,
College Park, MD 20742-4015}
\email{jda@math.umd.edu}
\thanks{The author was supported in part by NSF Grants \#DMS-0554278
  and \#DMS-0532393.}

\dedicatory{In memory of Fokko du Cloux}

\title[guide to the atlas software]
{Guide to the  Atlas Software:\\Computational Representation  Theory of Real Reductive Groups}

\maketitle

\sec{Introduction}
\label{i:introduction}

The Atlas of Lie Groups and Representations is a project in the
representation theory of real reductive groups.
The main goal of the {\tt atlas} computer software, currently under development, is
to compute the unitary dual of any real reductive Lie
group $G$. As a step in this direction it currently computes the
admissible representations of  $G$. 

The underlying mathematics of the software is described in 
\cite{algorithms}. 
See Sections 1 and 2 of \cite{algorithms}
for an overview of the algorithm.
This paper is a complement, and consists of a guide to the software
illustrated by numerous examples.

The software is currently in an early stage of development (version 0.3
as of April 2008). It is available from the Atlas web site {\tt
www.liegroups.org}. 

The {\tt help} command in the {\tt atlas} software is a another useful
source of information.  Also see the {\tt software} section of {\tt
www.liegroups.org}, including {\tt examples} and {\tt Tables of
Structure and Representation Theory}.  We plan to publish a manual for
the software at the time version 1.0 is released.

Marc van Leeuwen is currently working on the {\it
interpreter}, which provides simpler and more powerful input and output
methods for the software.  This is included in the current software
distribution; see the {\tt makefile} and the {\tt source/interpreter}
directory. We do not discuss this here since it is in a state of flux.
The input methods described in this paper will be available
indefinitely. The interpreter will provide additional output routines,
but there should only be minor changes to the output described
here.

The {\tt atlas} software was written by Fokko du Cloux.  Alfred Noel has
also contributed to it. In addition to writing the interpreter, Marc
van Leeuwen has made substantial changes to the software, and since
version 0.3 has been in charge of software development.

This paper grew out of lectures delivered at a conference in honor of Dragan
Mili\v{c}i\'{c} and Bill Casselman in July 2006. Bill had a
catalytic effect on the Atlas project, which got started at a
conference he organized in Montreal in 2002.

\subsec{How to read this paper}

You should start by downloading the {\tt atlas} software from the
Atlas web site {\tt www.liegroups.org} and installing it. Currently it
runs under unix, including Solaris and linux, Mac OSX, and Windows.
You should also have 
\cite{algorithms} handy.

A good approach to learning the software and \cite{algorithms}
is go through this paper, and do the examples using the
software. Each section begins with a brief summary of the relevant
material from \cite{algorithms}. We being very simply, and build up to
more complicated examples by the end. We assume some familiarity with
the theory of real reductive groups, for example see Chapters 1-8 of
\cite{overview}.

\sec{Defining basic data}
\label{s:basic}

The starting point of the algorithm is a complex reductive algebraic group $G$
together with an {\it inner class of real forms} of $G$.
The latter is determined by an involution $\gamma\in\Out(G)$, the group
of outer automorphisms of $G$. We refer to $(G,\gamma)$ as {\it basic
data} \cite[Section 4]{algorithms}.

We define a semidirect product $\G^\Gamma=G\rtimes\Gamma$, where
$\Gamma=\text{Gal}(\C/\R)=\{1,\sigma\}$.  
We let $\sigma$ act on $G$
by a ``distinguished'' involution $\tau\in\Aut(G)$ mapping to
$\gamma\in \Out(G)$, via the map $\Aut(G)\rightarrow \Out(G)$.  (A
distinguished involution is one which fixes a splitting datum
$(B,H,\{X_\alpha\})$; it is the ``most compact'' involution in this
inner class.) 
Write $G^\Gamma=\langle G,\delta\rangle$ where
$\delta=1\times\sigma$. 
See \cite[Section 5]{algorithms}.

A {\it real form} of $G$ in this inner class is a conjugacy class of
involutions $\theta\in \Aut(G)$ mapping to $\gamma\in\Out(G)$.
If $\theta$ is an involution of $G$ let $\sigma$ be an antiholomorphic
involution of $G$ commuting with $\theta$. Then $G(\R)=G^\sigma$ is a real
group in the usual sense, and $G(\R)^\theta$ is a maximal compact
subgroup of $G(\R)$. See \cite[Section 5]{algorithms}.

The reader should keep in mind the basic example of $\gamma=1$, so
$\GGamma=G\times\Gamma$. This is known as the {\it compact} or  {\it
  equal rank}  inner class; it contains the compact real form of $G$, and all
groups in this inner class contain a compact Cartan subgroup.

The first step in using the {\tt atlas} software is to define basic
data $(G,\gamma)$.
This proceeds in three steps:
\begin{enumerate}
\item Define a complex reductive Lie algebra $\g$, and let $G^*$ be
  the product of a complex torus and a simply connected, semisimple
  complex group with Lie algebra $\g$;
\item Choose a finite subgroup $A$ of $Z(G^*)$, and set $G=G^*/A$;
\item Choose an inner class of real forms of $G$.
\end{enumerate}

All three steps are accomplished by the {\tt type} command. We break
this up into three steps.

\subsec{Defining $\g$ and $G^*$}
\label{s:G^*}

A complex reductive Lie algebra is given by a list of types
{\tt A\_n,B\_n,\dots, E\_8, T\_n}, where {\tt T\_n} is the abelian Lie algebra $\C^n$.
In response to the  {\tt type} command, the software asks for the {\tt
Lie type:}. The user then enters such a list, with terms separated by 
a period. The order is irrelevant
here, although it plays a role in steps (2) and (3).
The entry {\tt T\_2} is the same as {\tt T\_1.T\_1}.

This defines the complex reductive Lie algebra $\g$, and group 
$G^*$. Here are some simple examples.

\begin{exampleplain}
We start with  $SL(2,\C)$:
\begin{verbatim}
main: type
Lie type: A1
\end{verbatim}
Here is $SL(2,\C)\times SL(2,\C)$:
\begin{verbatim}
main: type
Lie type: A1.A1
\end{verbatim}
Here is $SL(2,\C)\times \C^\times$:
\begin{verbatim}
main: type
Lie type: A1.T1
\end{verbatim}
and finally
$SL(2,\C)^2\times Spin(5,\C)\times Sp(4,\C)\times Sp(6,\C)\times \C^{\times4}$:
\begin{verbatim}
main: type
Lie type: A1.T1.B2.C2.T1.C3.T2.A1
\end{verbatim}
\end{exampleplain}

\subsec{Defining a complex group $G$}
\label{s:center}

The second step is to pick a finite subgroup of the center  $Z^*$
of $G^*$.

The center of each simple factor of $G$ is a finite cyclic group except in
type $D_{2n}$ in which case it is $\Ztoo\times\Ztoo$.
The finite cyclic group of order $n$ is denoted {\tt Z/n}, and is viewed
as the group
$$
\frac1n\Z/\Z
=\{\frac 0n,\frac 1n,\dots, \frac{n-1}n\}.
$$
The elements of finite order in $\C^\times$ are isomorphic to $\Q/\Z$,
and an element of this group is given by an element of $\Q$.

Once the user has given the Lie type, the software prompts the user
for a finite subgroup of $Z^*$, generated by  a set of elements of
$Z^*$. Each such element is given by a list of fractions, one for each
term in the Lie type (two for $D_{2n}$), separated by commas. Each element is given on a
single line; the empty line terminates this aspect of the input.
For example simply typing {\tt return} in response to this prompt
takes $A=1$ and $G=G^*$. Typing {\tt sc} has the same
effect. Typing {\tt ad} gives the adjoint group (actually a torus
times the adjoint group of the derived group of $G^*$).

\begin{exampleplain}
Here is the group $SL(2,\C)$:
\begin{verbatim}
main: type
Lie type: A1
elements of finite order in the center of the simply connected group:
Z/2
enter kernel generators, one per line
(ad for adjoint, ? to abort):
sc
enter inner class(es):
\end{verbatim}
Entering a carriage return alone also gives $SL(2,\C)$.
\end{exampleplain}

\begin{exampleplain}
To define  $PSL(2,\C)$ take $A$ to be the center of $SL(2,\C)$, which
is generated by the element of order $2$:
\begin{verbatim}
Lie type: A1
elements of finite order in the center of the simply connected group:
Z/2
enter kernel generators, one per line
(ad for adjoint, ? to abort):
1/2

enter inner class(es):
\end{verbatim}
Entering {\tt ad} instead of {\tt 1/2} has the same effect.
\end{exampleplain}

\begin{exampleplain}
  Here is $SO(10,\C)$, which is the quotient of $Spin(10,\C)$ by the
  element of $Z^*$ of order $2$:

\begin{verbatim}
main: type
Lie type: D5
elements of finite order in the center of the simply connected group:
Z/4
enter kernel generators, one per line
(ad for adjoint, ? to abort):
2/4

\end{verbatim}
\end{exampleplain}

\begin{exampleplain}
Often a reductive group is a quotient of $G^*$ by a diagonal
subgroup.
For example, here is $GL(2,\C)\simeq SL(2,\C)\times \C^*/\{(-I,-1)\}$:

\begin{verbatim}
main: type
Lie type: A1.T1
elements of finite order in the center of the simply connected group:
Z/2.Q/Z
enter kernel generators, one per line
(ad for adjoint, ? to abort):
1/2,1/2
\end{verbatim}
\end{exampleplain}

\begin{exampleplain}
\label{ex:so8C}

The center of $Spin(4n,\C)$ is not cyclic, and therefore requires two
terms.  For example, here is $SO(8,\C)\simeq Spin(8,\C)/A$ where $A$ is the diagonal
  subgroup of $Z^*=\Ztoo\times\Ztoo$:

\begin{verbatim}
main: type
Lie type: D4
elements of finite order in the center of the simply connected group:
Z/2.Z/2
enter kernel generators, one per line
(ad for adjoint, ? to abort):
1/2,1/2
\end{verbatim}
\end{exampleplain}

\begin{exampleplain}
\label{ex:SO8nondiag}
  We can also take  the quotient of $Spin(8,\C)$ by a non-diagonal 
  $\Ztoo$ subgroup:

\begin{verbatim}
main: type
Lie type: D8
elements of finite order in the center of the simply connected group:
Z/2.Z/2
enter kernel generators, one per line
(ad for adjoint, ? to abort):
1/2,0/2
\end{verbatim}
This is not isomorphic to $SO(8,\C)$.
\end{exampleplain}

This defines the group $G=G^*/A$ and completes step 2.

\subsec{Defining an inner class of real forms}
\label{s:inner}

We next define an inner class of real forms.
Recall this is determined by 
an involution in $\Out(G)$. 

The trivial element of $\Out(G)$ corresponds to the inner class of
real forms containing a compact Cartan subgroup.
This is the {\it compact}
inner class, and is denoted {\tt c}, and also {\tt e} for {\it equal
  rank}.
In particular if $G$ has a torus factor  its real points are
$S^1\times\dots\times S^1$.

Another natural inner class is that of the split real form.
This class is denoted {\tt s}; for a torus 
the real points are
$\R^\times\times\dots\times\R^\times$.
In many cases the classes
{\tt c} and {\tt s} are the same, for example if $\Out(G)=1$.

Now suppose $G=G_1\times G_1$. The outer automorphism switching the two
factors corresponds to the inner class of the real form $G_1(\C)$,
viewed as a real group. This class is denoted {\tt C}; in the case of
a torus this gives $\C^\times\times\dots\times\C^\times$ (viewed as a real
group). 

If $G$ is simple, then except in  type $D_{2n}$ the classes {\tt
c,s} exhaust every inner class. In type $D_{2n}$ the classes
{\tt c} and {\tt s} coincide, and there is another inner class 
denoted {\tt u}, for {\it unequal rank}. See example
\ref{ex:SO12_1}. This class is allowed in any
type for which the Dynkin diagram has a non-trivial automorphism. 

Now suppose $G=G^*$. An inner class of real forms of $G$ is specified
by choosing {\tt c,e,s,u} for each simple or $T^1$ factor, or {\tt C} for
each pair of identical (simple or $T^1$) factors, or factor of type 
$T^{2n}$.
In general an inner class of real forms of $G$ is  given by an {\it allowed} inner
class of real forms of $G^*$: the 
involutions in this inner class must factor to $G$. 

\bigskip
To summarize, to specify an inner class of real forms of $G$, give
a list of choices for each simple or torus factor, or pair of
identical entries in the case of $C$:

\begin{itemize}
\item {\tt c}: compact
\item {\tt e}: equal rank (same as {\tt c})
\item {\tt s}: split
\item {\tt u}: unequal rank
\item {\tt C}: complex (for an identical pair of entries).
\end{itemize}
The order of the choices corresponds to the order in which the simple
and torus factors of $G^*$ were specified.
If $G\ne G^*$ some choices may not be allowed.
\bigskip

Here are some examples.
Getting slightly ahead of ourself, the {\tt showrealforms} command lists
the real forms in the given inner class.
In the terminology of {\tt atlas} these are {\it weak} real forms; see
Sections \ref{s:real} and \ref{s:X}.
Also note that the first three inputs to the {\tt type} command can be
entered on a single line, provided  the second
input is {\tt ad} or {\tt sc}.

\begin{exampleplain}
The group $SL(2,\C)$ has two real forms $SL(2,\R)$ and $SU(2)$, both
of which are in the same inner class. Thus {\tt c=e=s} in this case.
\begin{verbatim}
empty: showrealforms
Lie type: A1 sc s
(weak) real forms are:
0: su(2)
1: sl(2,R)
\end{verbatim}
\end{exampleplain}
This example also illustrates a general principle of the software:
entering a command such as {\tt showrealforms} which requires other
input is allowed; the system will prompt the user for the
missing information.

\begin{exampleplain}
The smallest simple group with two inner classes is type $A_2$.
Here is the split inner class, with only one real form:

\begin{verbatim}
empty: showrealforms
Lie type: A2 sc s
(weak) real forms are:
0: sl(3,R)
\end{verbatim}
Here is the compact inner class, with two real forms:
\begin{verbatim}
empty: showrealforms
Lie type: A2 sc c
(weak) real forms are:
0: su(3)
1: su(2,1)
\end{verbatim}
\end{exampleplain}

\begin{exampleplain}
Here is $SL(3,\C)$ viewed as a real group:
\begin{verbatim}
empty: showrealforms
Lie type: A2.A2 sc C
(weak) real forms are:
0: sl(3,C)
\end{verbatim}
\end{exampleplain}

\begin{exampleplain}
If the group is a product specify the inner class as a list. 
For example here is the inner class of $SU(3)\times SL(3,\R)$: 
\begin{verbatim}
empty: showrealforms
Lie type: A2.A2 sc cs
(weak) real forms are:
0: su(3).sl(3,R)
1: su(2,1).sl(3,R)
\end{verbatim}
\end{exampleplain}
\noindent (The {\tt sc} means simply connected, and {\tt cs} means
compact$\times$split). 

\begin{exampleplain}
  Here is the inner class of $GL(2,\R)$:
\begin{verbatim}
empty: showrealforms
Lie type: A1.T1
elements of finite order in the center of the simply connected group:
Z/2.Q/Z
enter kernel generators, one per line
(ad for adjoint, ? to abort):
1/2,1/2

enter inner class(es): ss
(weak) real forms are:
0: su(2).gl(1,R)
1: sl(2,R).gl(1,R)
\end{verbatim}
\end{exampleplain}
\noindent Group 1 is $GL(2,\R)$. (Group 0 is the multiplicative group of the
quaternions.) 

\begin{exampleplain}
On the other hand here is the inner class of $U(1,1)$:
\begin{verbatim}
main: showrealforms
Lie type: A1.T1
elements of finite order in the center of the simply connected group:
Z/2.Q/Z
enter kernel generators, one per line
(ad for adjoint, ? to abort):
1/2,1/2

enter inner class(es): cc
(weak) real forms are:
0: su(2).u(1)
1: sl(2,R).u(1)
\end{verbatim}
\end{exampleplain}
\noindent Groups 0 and 1 are $U(2)$ and $U(1,1)$, respectively.

\begin{exampleplain}
\label{ex:SO12_1}
The group $SO(12,\C)$ has the following real forms, in two inner
classes. 
One inner class consists of the groups $SO(12,0)$, $S(10,2)$,
$SO(8,4)$, $SO(6,6)$ and
$SO^*(12)$. This is both the compact inner class (it contains $SO(12,0)$)
and split (it contains $SO(6,6)$), and is therefore the class {\tt
c=e=s}.
There is another inner class consisting of $SO(11,1), SO(9,3)$ and
$SO(7,5)$. This is the inner class {\tt u} of unequal rank.
This is an example where {\tt u} is needed.

Here is the compact and split inner class:
\begin{verbatim}
empty: showrealforms
Lie type: D6
elements of finite order in the center of the simply connected group:
Z/2.Z/2
enter kernel generators, one per line
(ad for adjoint, ? to abort):
1/2,1/2

enter inner class(es): e
(weak) real forms are:
0: so(12)
1: so(10,2)
2: so*(12)[1,0]
3: so*(12)[0,1]
4: so(8,4)
5: so(6,6)
\end{verbatim}
For an explanation of the two version of {\tt so*(12)} see
Example \ref{ex:SO12_2}.

The same input with {\tt u} in place of {\tt e} gives the
unequal rank inner class:
\begin{verbatim}
enter inner class(es): u
(weak) real forms are:
0: so(11,1)
1: so(9,3)
2: so(7,5)
\end{verbatim}
\end{exampleplain}

\begin{exampleplain}
Here is an example in which the inner class for $G^*$ is not defined
for $G$:

\begin{verbatim}
main: type
Lie type: A1.A1
elements of finite order in the center of the simply connected group:
Z/2.Z/2
enter kernel generators, one per line
(ad for adjoint, ? to abort):
1/2,0/2

enter inner class(es): C
sorry, that inner class is not compatible with the weight lattice
\end{verbatim}
What this means is: the automorphism which switches the two factors in $G^*=SL(2,C)\times
SL(2,\C)$ does not preserve $A$, and 
so does not factor to $G=PSL(2,\C)\times SL(2,\C)$.
\end{exampleplain}

\begin{exampleplain}
A more subtle example is the non-diagonal quotient
of $Spin(8,\C)$ of Example 
\ref{ex:SO8nondiag}:

\begin{verbatim}
main: type
Lie type: D4
elements of finite order in the center of the simply connected group:
Z/2.Z/2
enter kernel generators, one per line
(ad for adjoint, ? to abort):
1/2,0/2

enter inner class(es): u
sorry, that inner class is not compatible with the weight lattice
\end{verbatim}
\end{exampleplain}
\sec{Defining a real group}
\label{s:real}

We suppose the user has completed the {\tt type} command, which
defines $G$ and an inner class of real forms. We next specify a
particular real form of $G$ in the inner class.  
We discuss strong real forms in the next section.

The command {\tt showrealforms}, discussed in the preceding section, 
gives a list of the real forms of $\g$ in the given
inner class; the command {\tt realform} gives the same list, and the
user can choose the real form from the list.

The real forms of the classical Lie algebras  are given in the usual
notation. Some examples are {\tt sl(4,R)}, 
{\tt su(3,1)}, 
{\tt su(4)}, 
{\tt sl(4,H)}, 
{\tt sp(3,2)}, 
{\tt so(3,2)}, and
{\tt so*(10)}.
For a torus the real form is specified by {\tt gl(1,R)}, {\tt u(1)} or
{\tt gl(1,C)}.

For each exceptional group the real form is specified by specifying
the type of the maximal compact subgroup, except that the split real
form is denoted {\tt R}. For example the real forms of 
 $E_7$ are {\tt e7},
{\tt e7(e6.u(1))},
{\tt e7(so(12).su(2))}
and {\tt e7(R)}; otherwise known as {\it compact, Hermitian,
  quaternionic} and {\it split}, respectively.

\begin{exampleplain}
Here is the group $SL(2,\R)$:
\begin{verbatim}
empty: realform
Lie type: A1 sc s
(weak) real forms are:
0: su(2)
1: sl(2,R)
enter your choice: 1
\end{verbatim}
\end{exampleplain}

\begin{exampleplain}
Here is $PSL(2,\R)\simeq PGL(2,\R)\simeq SO(2,1)$.

\begin{verbatim}
main: type
Lie type: A1 ad s
main: realform
(weak) real forms are:
0: su(2)
1: sl(2,R)
enter your choice: 1
\end{verbatim}
Note that the real forms listed don't depend on the fact that the
complex group is $PSL(2,\C)$ here instead of $SL(2,\C)$; the real form
is determined by a real form of $\g$.
\end{exampleplain}

\begin{exampleplain}
\label{ex:SO12_2}
The equal rank forms of $SO(12,\C)$ were given in Example \ref{ex:SO12_1}:

\begin{verbatim}
0: so(12)
1: so(10,2)
2: so*(12)[1,0]
3: so*(12)[0,1]
4: so(8,4)
5: so(6,6)
\end{verbatim}
\end{exampleplain}
The two copies of {\tt so*(12)} illustrate
a technical point.
We say two real forms are equivalent if they are conjugate by $G$. See
\cite[Section 5]{algorithms}.
In the literature 
equivalence  is defined to be conjugacy by
$\Aut(G)$. If $G$ is simple the two notions agree in most
cases.

There are two real forms (in our sense) of $SO(12,\C)$ corresponding
to the real group $SO^*(12)$. In other words there are two isomorphic
subgroups of  $SO(12,\C)$, which are
not conjugate to each other. They {\it are} related by an outer
automorphism of $SO(12,\C)$.
(In the Kac classification of real forms each is labelled by a $2$ on
one of the branches of the fork in the Dynkin diagram).
These two groups are denoted $SO^*(12)[1,0]$ 
and $SO^*(12)[0,1]$. 
The fact that there are two copies of $SO^*(12)$ manifests itself in
the 
representation theory of $SO(p,q)$. See Example \ref{ex:so12}.

In fact we can even see this distinction in terms of the structure
theory of the groups themselves.
If $G=Spin(12,\C), SO(12,\C)$ or $PSO(12,\C)$ the two real groups
locally isomorphic to $SO^*(12)$ are
isomorphic, and interchanged by an outer automorphism of $G$. 
However  in the non-diagonal quotient of  $Spin(12,\C)$ this
fails, as is illustrated in the next example.

\begin{exampleplain}
Here are two versions of $SO^*(12)$ in the non-diagonal quotient of $Spin(12,\C)$:
\begin{verbatim}
main: type
Lie type: D6
elements of finite order in the center of the simply connected group:
Z/2.Z/2
enter kernel generators, one per line
(ad for adjoint, ? to abort):
1/2,0/2

enter inner class(es): c
main: realform
(weak) real forms are:
0: so(12)
1: so(10,2)
2: so*(12)[1,0]
3: so*(12)[0,1]
4: so(8,4)
5: so(6,6)
enter your choice: 2
real: components
component group is (Z/2)^1
real: realform
(weak) real forms are:
0: so(12)
1: so(10,2)
2: so*(12)[1,0]
3: so*(12)[0,1]
4: so(8,4)
5: so(6,6)
enter your choice: 3
real: components
group is connected
\end{verbatim}
Note that these two groups are {\it not} isomorphic.
This illustrates the command {\tt components}, which
gives the component group of the real group (an elementary abelian two-group).
\end{exampleplain}

In $D_4$ the situation is even more interesting. Let $G=Spin(8,\C)$
and consider the  three real forms,  $Spin(6,2)$ and the two
versions of $Spin^*(8)$. In fact these are isomorphic;
$Spin(6,2)\simeq Spin^*(8)$, and these three groups are interchanged
by the outer automorphism group of $Spin(8,\C)$, which is $S_3$.
Similar statements hold in $PSO(8,\C)$. We leave it to the reader to
analyze real forms of $SO(8,\C)$ and the non-diagonal quotient of
$Spin(8,\C)$, in which this symmetry is broken in different ways.
See Examples \ref{ex:so8C} and  \ref{ex:SO8nondiag}.

In general in $SO(4n,\C)$ (type $D_{2n}$) there are two non-conjugate copies of
$SO^*(4n)$. However in $SO(4n+2,\C)$  (type $D_{2n+1}$) 
all copies of $SO^*(4n+2)$ are in fact conjugate by $SO(2n+2,\C)$. 
This is illustrated by the next example.

\begin{exampleplain}
\label{ex:so10}
Here are the equal rank  real forms of $SO(10,\C)$:

\begin{verbatim}
empty: showrealforms
Lie type: D5
elements of finite order in the center of the simply connected group:
Z/4
enter kernel generators, one per line
(ad for adjoint, ? to abort):
2/4

enter inner class(es): c
(weak) real forms are:
0: so(10)
1: so(8,2)
2: so*(10)
3: so(6,4)
\end{verbatim}
\end{exampleplain}
Note that there is only one real form $SO^*(10)$.

\sec{The Spaces $\X$ and $\X^r$}
\label{s:X}

Fix basic data $(G,\gamma)$ as in Section \ref{s:basic}.

By a {\it strong involution} for $(G,\gamma)$ we mean an element $\xi\in
G^\Gamma\bs G$ satisfying $\xi^2\in Z$ ($Z=Z(G)$ is the center of $G$).
Let $\I$ be the set of strong involutions.
A {\it strong real form} is an equivalence class of strong
involutions (equivalence is by conjugacy by $G$).
If $\xi$ is a strong involution then 
$\theta_\xi=\int(\xi)$
is an involution of $G$ in the inner class of $\gamma$,
and this gives  a surjective map 
from strong real forms  to real forms in this inner class. 
If
$G$ is adjoint this map is bijective, but not otherwise.
Let $K_\xi=G^{\theta_\xi}$; this is the complexified maximal compact
subgroup of the corresponding real group $G(\R)$.
See Section \ref{s:basic}. 

Recall (Section \ref{s:basic}) $\delta=1\times\sigma$ is the
distinguished element of $G^\Gamma$.
We fix once and for all a Cartan subgroup $H$ of $G$ stable by
$\int(\delta)$. 
Define
\begin{equation}
\label{e:X}
\begin{aligned}
\tX&=\{g\in \Norm_{G^\Gamma\bs G}(H)\,|\, g^2\in Z\}\\
\X&=\tX/H
\end{aligned}
\end{equation}
(the quotient is by the conjugation action of $H$).
See \cite[Section 9]{algorithms}.
The space $\X$ is our main combinatorial object, 
known as the {\it (one sided) parameter space}.

Mathematically the space $\X$ is a natural object. From the point of
view of computations it has the disadvantage that (if $G$ is not semisimple) it may be infinite,
and the fibers of the map from strong real forms to real forms may be
infinite. The software works exclusively with the {\it reduced
  parameter space} $\X^r$ which we now describe. See \cite[Section
13]{algorithms}.  

 Let $\theta=\int(\delta)$ and 
choose a set 
$Z^r\subset Z$ of representatives of
\begin{equation}
\{z\in Z\,|\,\theta(z)=z\}/\{z\theta(z)\,|\,z\in Z\}.
\end{equation}
Define reduced versions of $\I,\tX$ and $\X$:
\begin{equation}
\begin{aligned}
\I^r&=\{g\in\GGamma\bs G\,|\, g^2\in Z^r\}\\
\tX^r&=\{g\in\Norm_{\GGamma\bs G}(H)\,|\, g^2\in Z^r\}, \X^r=\tX^r/H.
\end{aligned}
\end{equation}
Henceforth we define a strong real form in this revised sense: $\I^r$
is the space of strong real forms, 
with equivalence given by conjugacy by $G$ as before.
The map from strong real forms 
to real forms is still surjective. The
fibers of this map are finite, and 
$\Xr$ is a finite set.

Choose a set $\{\xi_i\,|\,1\le i\le n\}\subset \tX$ of representatives of strong real forms
(in the sense of the reduced parameter space). That is
\begin{equation}
\label{e:I}
\{\xi_1,\dots, \xi_n\}
\overset{1-1}\longleftrightarrow
\I^r/G.
\end{equation}
For each $i$ let $x_i$ be the image of
$\xi_i$ in 
$\Xr$,
$\theta_i=\int(\xi_i)$ and
$K_i=G^{\theta_i}$. See \cite[5.15]{algorithms}.

Let $W^\Gamma=\Norm_{\GGamma}(H)/H$ and 
$$
\I_W=\{\tau\in W^\Gamma\bs W\,|\, \tau^2=1\}. 
$$
We continue to write $\delta$ for the image of $\delta\in \GGamma$ in
$W^\Gamma$, and $\theta$ for the involution $\int(\delta)$ of $W$.
Then $W^\Gamma=\langle W,\delta\rangle$.
By a {\it twisted involution} in $W$ we mean an element $w\in W$
satisfying $w\theta(w)=1$, and we say two such elements $w,w'$ are {\it
  twisted}-conjugate if $w'=yw\theta(y\inv)$ for some $y\in W$.
The map $w\rightarrow w\delta$ is a bijection between the space of
twisted involutions and $\I_W$, taking twisted-conjugacy to ordinary
conjugacy by $W$. We refer to $\I_W$ as the space of twisted
involutions, and pass back and forth between the two notions.
See \cite[9.14]{algorithms}.  

There is a natural
surjective map $p:\X^r\rightarrow\I_W$ \cite[Lemma
9.12]{algorithms}.  The fiber  over an element
$\tau\in \I_W$ is denoted $\Xr_\tau$.

The group $\Norm_{G}(H)$ acts by conjugation on $\tXr$, and
this factors to an action of $W$ on $\Xr$.
Because of its
relationship to the cross action of 
\cite[Definition 8.3.1]{green}
we call this the {\it cross action} of $W$, and denote it $w\times x$.
This action is equivariant for $p$ and the conjugation action of $W$
on $\I_W$: $p(w\times x)=wp(x)w\inv$.

\sec{Cartan subgroups}
\label{s:cartan}

Fix basic data $(G,\gamma)$.
Let $\theta_{qs}$ be a quasisplit involution in this inner class
and let $K_{qs}=G^{\theta_{qs}}$. Thus $K_{qs}$ is the complexified maximal
  compact subgroup of the quasisplit form $G_{qs}(\R)$ of $G$. 
The conjugacy classes of
Cartan subgroups  of $G_{qs}(\R)$, equivalently the $K_{qs}$-conjugacy
classes of $\theta_{qs}$-stable Cartan subgroups of $G$, are in
natural bijection with 
$\I_W/W$, conjugacy classes of twisted involutions in $W$.
The conjugacy classes of  Cartan subgroups for any real form of $G$
are a subset of those for the quasisplit real form.
See \cite[Propositions 12.9 and 12.12]{algorithms}.

The {\tt cartan} command gives a list of Cartan subgroups for a given
inner form. For each Cartan subgroup it gives the following
information. See the {\tt help} file for the {\tt cartan} command for
more details.

First it gives the structure of  the Cartan subgroup 
as a real torus: $H(\R)$ is isomorphic to 
$(\R^\times)^a\times
(S^1)^b\times (\C^\times)^c$ for integers $(a,b,c)$; these are the
{\tt split}, {\tt compact} and {\tt complex} entries in the output of
{\tt cartan}.

Each conjugacy class of twisted involutions contains a unique
canonical representative $\tau=w\delta$. The next line of output is a
reduced expression for $w$.  

The first entry on the next line is the number of twisted involutions
in this conjugacy class. For an explanation of {\tt fiber rank} 
and {\tt \#X\_r} see Section \ref{s:geometryofX}.

Associated to any Cartan subgroup are the sets of real and imaginary
roots (see Section \ref{s:kgb}),  each of which is a root system. These are given by {\tt
  imaginary root system} and {\tt real root system}, respectively. The {\tt complex
  factor} line in the output has to do with the Weyl group; see
Section \ref{s:weyl}. 

The lines beginning {\tt real form\dots} give
information about $\X^r_\tau$. For now we
observe only that for a given Cartan subgroup only the real forms which contain this
Cartan subgroup are displayed.
See Section \ref{s:adjoint}.

\begin{exampleplain}
Let $G=SL(2,\C)$, which has a unique inner class of real forms.
Then $\I_W=W\times\Gamma$; we can drop $\Gamma$ and write
$\I_W/W=\{1,s\}$. There are two conjugacy classes of Cartan subgroups
of $SL(2,\R)$, compact and split.
\begin{verbatim}
empty: cartan
Lie type: A1 sc s
(weak) real forms are:
0: su(2)
1: sl(2,R)
enter your choice: 1
Name an output file (return for stdout, ? to abandon):

Cartan #0:
split: 0; compact: 1; complex: 0
canonical twisted involution:
twisted involution orbit size: 1;  fiber rank: 1;  #X_r: 2
imaginary root system: A1
real root system is empty
complex factor is empty
real form #1: [0] (1)
real form #0: [1] (1)


Cartan #1:
split: 1; compact: 0; complex: 0
canonical twisted involution: 1
twisted involution orbit size: 1;  fiber rank: 0;  #X_r: 1
imaginary root system is empty
real root system: A1
complex factor is empty
real form #1: [0] (1)
\end{verbatim}

{\tt Cartan \#0} is always the fundamental (most compact) Cartan
subgroup, in this case $S^1$. The corresponding twisted involution is
the identity, and its orbit is itself. Both real forms contain this
Cartan subgroup. 

{\tt Cartan \#1} is the split Cartan subgroup $\R^\times$, and this
occurs only in real form {\tt \#1}, i.e. $SL(2,\R)$.

\end{exampleplain}

\begin{exampleplain}
\label{ex:sl2C}

Here is the complex group $SL(2,\C)$:

\begin{verbatim}
empty: cartan
Lie type: A1.A1 sc C
there is a unique real form: sl(2,C)
Name an output file (return for stdout, ? to abandon):

Cartan #0:
split: 0; compact: 0; complex: 1
canonical twisted involution:
twisted involution orbit size: 2;  fiber rank: 0;  #X_r: 2
imaginary root system is empty
real root system is empty
complex factor: A1
real form #0: [0] (1)
\end{verbatim}
There is one real form, and one  conjugacy class of Cartan subgroups, 
isomorphic to  $\C^\times$.
For a complex group the twisted involutions are in bijection with the
Weyl group, so in this case there are $2$.
\end{exampleplain}

\begin{exampleplain}
Here are the Cartan subgroups of $Sp(4,\R)$. See \cite[Example
14.19]{algorithms}. 

\begin{verbatim}
empty: cartan
Lie type: C2 sc s
(weak) real forms are:
0: sp(2)
1: sp(1,1)
2: sp(4,R)
enter your choice: 2
Name an output file (return for stdout, ? to abandon):

Cartan #0:
split: 0; compact: 2; complex: 0
canonical twisted involution:
twisted involution orbit size: 1;  fiber rank: 2;  #X_r: 4
imaginary root system: B2
real root system is empty
complex factor is empty
real form #2: [0,1] (2)
real form #1: [2] (1)
real form #0: [3] (1)

Cartan #1:
split: 0; compact: 0; complex: 1
canonical twisted involution: 2,1,2
twisted involution orbit size: 2;  fiber rank: 0;  #X_r: 2
imaginary root system: A1
real root system: A1
complex factor is empty
real form #2: [0] (1)
real form #1: [1] (1)

Cartan #2:
split: 1; compact: 1; complex: 0
canonical twisted involution: 1,2,1
twisted involution orbit size: 2;  fiber rank: 1;  #X_r: 4
imaginary root system: A1
real root system: A1
complex factor is empty
real form #2: [0] (1)

Cartan #3:
split: 2; compact: 0; complex: 0
canonical twisted involution: 1,2,1,2
twisted involution orbit size: 1;  fiber rank: 0;  #X_r: 1
imaginary root system is empty
real root system: B2
complex factor is empty
real form #2: [0] (1)
\end{verbatim}
\end{exampleplain}

There are four conjugacy classes of Cartan subgroups, isomorphic to
$S^1\times S^1, \C^\times,S^1\times\R^\times$ and $\R^\times\times\R^\times$.
All four are contained in the split group $Sp(4,\R)$; two of them are
contained in $Sp(1,1)$, and only the compact Cartan subgroup occurs in
$Sp(2,0)$.

\begin{exampleplain}
\label{ex:E6f4}
Here is the real form of $E_6$ with $K$ of type $F_4$:

\begin{verbatim}
empty: cartan
Lie type: E6 sc s
(weak) real forms are:
0: e6(f4)
1: e6(R)
enter your choice: 0
Name an output file (return for stdout, ? to abandon):

Cartan #0:
split: 0; compact: 2; complex: 2
canonical twisted involution:
twisted involution orbit size: 45;  fiber rank: 2;  #X_r: 180
imaginary root system: D4
real root system is empty
complex factor: A2
real form #1: [0,1,2] (3)
real form #0: [3] (1)
\end{verbatim}
This group has a unique conjugacy class of Cartan subgroups.
See Example \ref{ex:e6}.
\end{exampleplain}

We discuss Cartan subgroups of the Classical groups.

\begin{exampleplain}
Type $A_{n-1}$.
It is most convenient to take $G=GL(n,\C)$.

First let $\gamma=1$, so $W^\Gamma=S_n\times\Gamma$. 
It is well known that the conjugacy classes of involutions in $S_n$
are parametrized by  ordered pairs
$(a,b)\in\mathbb N^2$  satisfying $a+2b=n$: in cycle notation take
$w=(1,2)(3,4)\dots(2b-1,2b)$.

The quasisplit group in this inner class is the  unitary
group $U(m,m)$ or $U(m+1,m)$.
The Cartan subgroup corresponding to $(a,b)$ is isomorphic to
$(S^1)^a\times (\C^\times)^b$; the identity element
corresponds to the compact Cartan subgroup.

Now suppose $\gamma$ is given by the unique non-trivial automorphism of the
Dynkin diagram ($n\ge 3$).
The split real form is $GL(n,\R)$.
It turns out that the twisted involutions in $W$
are also parametrized by pairs $(a,b)$ with $a+2b=n$; in this case the
corresponding Cartan subgroup is isomorphic to
$(\R^\times)^a\times(\C^\times)^b$.
(This is an aspect of Vogan duality; see Section \ref{s:rep} and
\cite[1.35 and Corollary 10.9]{algorithms}.)
\end{exampleplain}

\begin{exampleplain}
Types $B_n$ and $C_n$.

In this   case $\gamma$ is necessarily trivial, and 
the Cartan subgroups of 
$SO(n+1,n)$ or $Sp(2n,\R)$ are parametrized by conjugacy classes of
involutions in $W\simeq S_n\ltimes\Ztoo^n$.
These 
are parametrized by 
$(a,b,c)$ with $a+b+2c=n$, and the corresponding Cartan subgroup is
isomorphic to $(S^1)^a\times(\R^\times)^b\times(\C^\times)^c$. 

Is is interesting to consider the Hasse diagram of these Cartan
subgroups. This is the graph, with one node for each Cartan subgroup,
and an edge for each Cayley transform (cf.~Section \ref{s:kgb}) relating two Cartan subgroups. 
We make a node black if the 
corresponding Cartan
subgroup is the most split Cartan subgroup of a real form of
$G$.    The nodes corresponding to Cartan subgroups with the same
split rank are arranged in the same row, and each row is numbered with
the value of the split rank.  Such diagrams for $SO(7,6)$ and $Sp(12,\R)$ are given in Figure
\ref{fig:1}.

\newpage

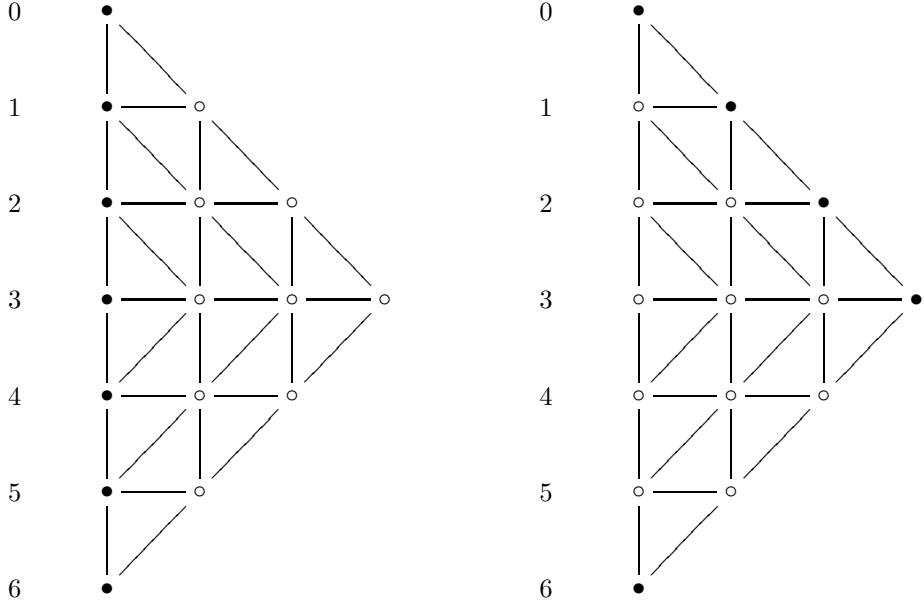
\begin{figure}[ht]
\[
\xymatrix{
0&{\bullet}\ar@{-}[d]\ar@{-}[dr]\\
1&{\bullet}\ar@{-}[d]\ar@{-}[dr]\ar@{-}[r]&{\circ}\ar@{-}[d]\ar@{-}[dr]\\
2&{\bullet}\ar@{-}[d]\ar@{-}[dr]\ar@{-}[r]&{\circ}\ar@{-}[d]\ar@{-}[dr]\ar@{-}[r]&{\circ}\ar@{-}[d]\ar@{-}[dr]\\
3&{\bullet}\ar@{-}[d]\ar@{-}[r]&{\circ}\ar@{-}[d]\ar@{-}[r]&{\circ}\ar@{-}[d]\ar@{-}[r]&{\circ}&\\
4&{\bullet}\ar@{-}[ur]\ar@{-}[d]\ar@{-}[r]&{\circ}\ar@{-}[ur]\ar@{-}[d]\ar@{-}[r]&{\circ}\ar@{-}[ur]\\
5&{\bullet}\ar@{-}[d]\ar@{-}[ur]\ar@{-}[r]&{\circ}\ar@{-}[ur]\\
6&{\bullet}\ar@{-}[ur]\\
}
\qquad
\xymatrix{
0&{\bullet}\ar@{-}[d]\ar@{-}[dr]\\
1&{\circ}\ar@{-}[d]\ar@{-}[dr]\ar@{-}[r]&{\bullet}\ar@{-}[d]\ar@{-}[dr]\\
2&{\circ}\ar@{-}[d]\ar@{-}[dr]\ar@{-}[r]&{\circ}\ar@{-}[d]\ar@{-}[dr]\ar@{-}[r]&{\bullet}\ar@{-}[d]\ar@{-}[dr]\\
3&{\circ}\ar@{-}[d]\ar@{-}[r]&{\circ}\ar@{-}[d]\ar@{-}[r]&{\circ}\ar@{-}[d]\ar@{-}[r]&{\bullet}&\\
4&{\circ}\ar@{-}[ur]\ar@{-}[d]\ar@{-}[r]&{\circ}\ar@{-}[ur]\ar@{-}[d]\ar@{-}[r]&{\circ}\ar@{-}[ur]\\
5&{\circ}\ar@{-}[d]\ar@{-}[ur]\ar@{-}[r]&{\circ}\ar@{-}[ur]\\
6&{\bullet}\ar@{-}[ur]
}
\]
\caption{Cartan subgroups for $SO(7,6)$ and $Sp(12,\R)$.}
\label{fig:1}
\end{figure}
\end{exampleplain}

\begin{exampleplain}
\label{ex:dn}
Type $D_n$.  

In this case (for $n\ge 2$) there are two choices of
$\gamma$, corresponding to the quasisplit groups $SO(n,n)$ and
$SO(n+1,n-1)$
($\gamma=1$ corresponds to $SO(n,n)$ if $n$ is even, and $SO(n+1,n-1)$
if $n$ is odd).
It is most convenient to group these two real forms together. 
Then the Cartan subgroups are parametrized by  $(a,b,c)$ with
$a+b+2c=n$, except that $(0,0,c)$ is counted twice.
(These two Cartan subgroups are conjugate by the outer automorphism of
$SO(n,n)$ coming from $O(n,n)$.)
Again the corresponding Cartan subgroup is 
isomorphic to $(S^1)^a\times(\R^\times)^b\times(\C^\times)^c$.
If $a$ is even this Cartan subgroup occurs in $SO(n,n)$ and in
$SO(n+1,n-1)$ otherwise.  Diagrams for $SO(6,6)$ and $SO(7,5)$
are given in Figure \ref{fig:2}.

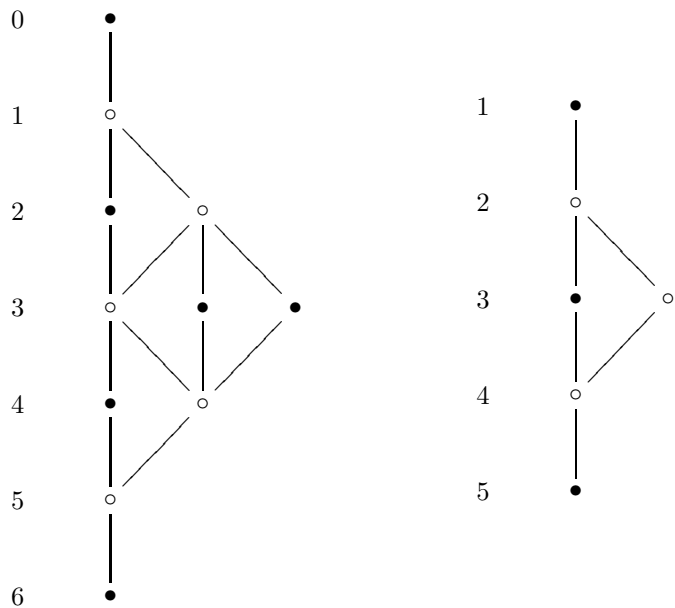
\begin{figure}[ht]
\begin{equation*}
\xymatrix{
0&{\bullet}\ar@{-}[d]\\
1&{\circ}\ar@{-}[d]\ar@{-}[dr]\\
2&{\bullet}\ar@{-}[d]&{\circ}\ar@{-}[d]\ar@{-}[dr]\\
3&{\circ}\ar@{-}[d]\ar@{-}[ur]\ar@{-}[dr]&{\bullet}\ar@{-}[d]&{\bullet}\\
4&{\bullet}\ar@{-}[d]&{\circ}\ar@{-}[ur]\\
5&{\circ}\ar@{-}[d]\ar@{-}[ur]\\
6&{\bullet}}
\qquad
\qquad
\qquad
\xymatrix{\\
1&{\bullet}\ar@{-}[d]\\
2&{\circ}\ar@{-}[d]\ar@{-}[dr]\\
3&{\bullet}\ar@{-}[d]&{\circ}\ar@{-}[dl]\\
4&{\circ}\ar@{-}[d]\\
5&{\bullet}}
\end{equation*}
\caption{Cartan subgroups for $SO(6,6)$ and $SO(7,5)$.}
\label{fig:2}
\end{figure}
\bigskip

\begin{exampleplain}
\label{ex:realtorus}
Note that the Cartan subgroups in $G(\R)$ only depend on $W^\Gamma$,
and are therefore independent of isogeny (this is not true in the
p-adic case). 
While  the list of Cartan subgroups is independent of
isogeny, the description as a real torus is not, and the nature of the
torus can change unexpectedly under isogenies.

For example suppose $G$ is of type $D_2\simeq A_1\times A_1$ and
$\gamma\ne 1$. There a unique conjugacy class of Cartan subgroups in
this case. 
In the simply connected case $G(\R)\simeq SL(2,\C)$ and
the Cartan subgroup is isomorphic to $\C^\times$ (cf.~Example
\ref{ex:sl2C}). 

For $SO(3,1)$ we get $S^1\times\R^\times$:
\begin{verbatim}
split: 1; compact: 1; complex: 0
\end{verbatim}
and  for $PSO(3,1)\simeq PSL(2,\C)$ we get $\C^\times$ again:
\begin{verbatim}
split: 0; compact: 0; complex: 1
\end{verbatim}
Note that $G(\R)$ is a connected complex group if $G(\C)=Spin(4,\C)$ or $PSO(4,\C)$, but not $SO(4,\C)$.
In fact we have:

\bigskip

\begin{tabular}{lll}
$G$&$G(\R)$& $H(\R)$\\
$Spin(4,\C)\simeq SL(2,\C)\times SL(2,\C)$ & $Spin(3,1)\simeq SL(2,\C)$ & $\C^\times$\\
$SO(4,\C)$ & $SO(3,1)$ & $\R^\times\times S^1$\\
$PSO(4,\C)\simeq PSL(2,\C)\times PSL(2,\C)$ & $PSO(3,1)\simeq PSL(2,\C)$ & $\C^\times$\\
\end{tabular}

\medskip

\end{exampleplain}

\begin{remarkplain}
If $G(\R)$ is a real form of either $GL(n,\C),SL(n,\C)$, $SO(n,\C)$ or $Sp(2n,\C)$
then two Cartan subgroups are isomorphic if and only if they are
conjugate by $G(\R)$, or by an automorphism of $G(\R)$ in the case of
$(\C^\times)^n\subset SO(2n,2n)$ (cf.~Example \ref{ex:dn}).
It is perhaps surprising that this fails badly
for isogenous groups, as the following
example shows.
\end{remarkplain}

\begin{exampleplain}
Let $G(\R)=Spin(n,n)$ or $PSO(n,n)$ with $n$ even ($n\ge 4$). There
are three non-conjugate Cartan subgroups isomorphic to $\R^\times\times S^1\times
(\C^\times)^{\frac n2-1}$. 
Two of these are interchanged by an outer automorphism of $G(\R)$.
(The corresponding Cartan subgroups of $SO(n,n)$ are isomorphic to
$(\C^\times)^{\frac n2}$.) The third one is not related to the others:
it has a different real Weyl group (see the  {\tt complex factor} entry).

Here are the three Cartan subgroups in question for $Spin(6,6)$.

\begin{verbatim}
empty: cartan
Lie type: D6 sc s
(weak) real forms are:
0: so(12)
1: so(10,2)
2: so*(12)[1,0]
3: so*(12)[0,1]
4: so(8,4)
5: so(6,6)
enter your choice: 5
Name an output file (return for stdout, ? to abandon):
...
Cartan #4:
split: 1; compact: 1; complex: 2
canonical twisted involution: 3,4,5,6,4,3,2,3,4,5,6,4,3,1,2,3,4,5,6,4,3,2,1
twisted involution orbit size: 180;  fiber rank: 1;  #X_r: 360
imaginary root system: A1.A1.A1
real root system: A1.A1.A1
complex factor: A1
real form #5: [0] (1)
real form #4: [1] (1)

Cartan #5:
split: 1; compact: 1; complex: 2
canonical twisted involution: 6,4,5,3,4,6,2,3,4,5,1,2,3,4,6
twisted involution orbit size: 60;  fiber rank: 1;  #X_r: 120
imaginary root system: A1.A1.A1
real root system: A1.A1.A1
complex factor: A2
real form #5: [0] (1)
real form #3: [1] (1)

Cartan #6:
split: 1; compact: 1; complex: 2
canonical twisted involution: 5,4,6,3,4,5,2,3,4,6,1,2,3,4,5
twisted involution orbit size: 60;  fiber rank: 1;  #X_r: 120
imaginary root system: A1.A1.A1
real root system: A1.A1.A1
complex factor: A2
real form #5: [0] (1)
real form #2: [1] (1)
...
\end{verbatim}

\end{exampleplain}
A similar phenomenon holds in the split real forms of $F_4$ and $G_2$.
\end{exampleplain}

\sec{Weyl Groups}
\label{s:weyl}

Fix basic data $(G,\gamma)$, $\xi\in\tX$, and set
$K=G^{\theta_\xi}$ (cf.~Section \ref{s:X}).
The ``real'' Weyl group $W(K,H)=\Norm_K(H)/H\cap K$ plays an important
role.  It is isomorphic to $W(G(\R),H(\R))=\Norm_{G(\R)}(H(\R))/H(\R)$
where $G(\R)$ is a real form of $G$ corresponding to $K$, and $H(\R)$
is the corresponding real form of $H$. See \cite[Section 12]{algorithms}.

We briefly recall some constructions from \cite[Section
8]{algorithms}, also see \cite[Proposition 4.16]{ic4}. 
Let $\tau$ be the image of $\xi$ in $\I_W$.
We have
\begin{equation}
W(K,H)\simeq (W_C)^\tau\ltimes (W(M\cap K,H)\times W_r).
\end{equation}
Here 
\begin{itemize}
\item 
$W_r$ is the Weyl group of the system of real roots;
\item $(W_C)^\tau$ is the Weyl group of a certain root system
constructed using complex roots (\cite[Proposition 3.12]{ic4});
\item $W(K\cap M,H)\simeq W_{i,c}\ltimes\mathcal A(H)\subset W_i$,
where  $W_i$ is the Weyl group of the roots system of imaginary roots,
$W_{i,c}$ is the Weyl group of the root system of compact
  imaginary roots, and $\mathcal A(H)$ is a certain elementary abelian
  two-group.   
\end{itemize}
See Section \ref{s:kgb} for the definition  of real and imaginary
roots.

To describe $W(K,H)$ it is therefore sufficient to describe $W_r,
W_{i,c}, (W_C)^\tau$ and $\mathcal A(H)$. The first three are Weyl
groups of root systems.
The {\tt realweyl} command describes these three root systems,
denoted {\tt  W\_r, W\_ic}, and {\tt W\^{}C}
 respectively; and $\mathcal A(H)$,
denoted {\tt A}.

\begin{exampleplain}
As usual we start with $SL(2,\R)$:
\begin{verbatim}
empty: realweyl
Lie type: A1 sc s
(weak) real forms are:
0: su(2)
1: sl(2,R)
enter your choice: 1
cartan class (one of 0,1): 0
Name an output file (return for stdout, ? to abandon):

real weyl group is W^C.((A.W_ic) x W^R), where:
W^C is trivial
A is trivial
W_ic is trivial
W^R is trivial
\end{verbatim}
This is the compact Cartan subgroup, for which 
$W(K,H)$ is trivial. Here is the split Cartan subgroup:

\begin{verbatim}
real: realweyl
cartan class (one of 0,1): 1
Name an output file (return for stdout, ? to abandon):

real weyl group is W^C.((A.W_ic) x W^R), where:
W^C is trivial
A is trivial
W_ic is trivial
W^R is a Weyl group of type A1

generators for W^R:
1
\end{verbatim}
For the split Cartan subgroup the only non-trivial factor is $W_r$, and
$W(K,H)=W_r\simeq\Ztoo$. 
  
\end{exampleplain}

\begin{exampleplain}
There is a small change when we compute the Weyl group
for the compact Cartan subgroup of $PSL(2,\R)$ instead:

\begin{verbatim}
empty: realweyl
Lie type: A1 ad s
(weak) real forms are:
0: su(2)
1: sl(2,R)
enter your choice: 1
cartan class (one of 0,1): 0
Name an output file (return for stdout, ? to abandon):

real weyl group is W^C.((A.W_ic) x W^R), where:
W^C is trivial
A is an elementary abelian 2-group of rank 1
W_ic is trivial
W^R is trivial

generators for A:
1
\end{verbatim}
In this case $\mathcal A(H)=\Ztoo$, so $W(K,H)=\Ztoo$. Recall $PSL(2,\R)\simeq
SO(2,1)$ is disconnected; the non-trivial Weyl group element is given
by an element in the non-identity component.
\end{exampleplain}

\begin{exampleplain}
Here is an example, for split groups of type $D_{2m}$,
in which the group $\mathcal A(H)$ is quite large.
We take the Cartan subgroup 
$(S^1)^2\times(\C^\times)^{m-1}$ (cf. Example \ref{ex:dn}).
In the notation at the beginning of \cite[Section 12]{algorithms} we have
\begin{equation}
\begin{aligned}
\Delta_i&=A_1^{m-1}\times D_2\simeq A_1^{m+1}\\
\Delta_r&=A_1^{m-1}\\
\Delta_C&=A_{m-2}\times A_{m-2}\\
W_i&\simeq \Ztoo^{m+1},\,W_r\simeq\Ztoo^{m-1}\\
W_{i,c}&=1\\
(W_C)^\tau&\simeq S_{m-1}\\
W^\tau&\simeq S_{m-1}\ltimes [\Ztoo^{m+1}\times\Ztoo^{m-1}].
\end{aligned}
\end{equation}
Here $S_{m-1}$ acts trivially on the final two factors of
$\Ztoo^{m+1}$, and $W^\tau$ is the centralizer of $\tau$ in $W$; this
group contains $W(K,H)$.
So far this discussion is independent of isogeny. 

Since $W_{i,c}=1$ the real Weyl group is the same as
$W^\tau$, with the factor $W_i=\Ztoo^{m+1}$ replaced by 
$\mathcal A(H)\subset W_i$.
The group $\mathcal A(H)$ depends on the isogeny. If $G$ is adjoint it
is equal to $W_i$. The smallest possible value of $\mathcal A(H)$
occurs when $G$ is simply connected.
It isn't as easy to compute $\mathcal A(H)$ in this case; the software
will tell us.

For example consider the split group $Spin(8,8)$ of type $D_8$ ($m=4$). The Cartan
class is {\tt \#5}:
\begin{verbatim}
Cartan #5:
split: 0; compact: 2; complex: 3
canonical twisted involution: 6,7,8,6,5,[truncated]
twisted involution orbit size: 3360;  fiber rank: 2;  #X_r: 13440
imaginary root system: A1.A1.A1.A1.A1
real root system: A1.A1.A1
complex factor: A2
...
\end{verbatim}
and here is the relevant output from {\tt realweyl}:
\begin{verbatim}
real weyl group is W^C.((A.W_ic) x W^R), where:
W^C is isomorphic to a Weyl group of type A2
A is an elementary abelian 2-group of rank 3
W_ic is trivial
W^R is a Weyl group of type A1.A1.A1
\end{verbatim}
Thus $\mathcal A(H)\simeq \Ztoo^{m-1}$.
\end{exampleplain}

If $H$ is a fundamental (most compact) Cartan subgroup then $H\cap K$ is a Cartan
subgroup of $K$, and 
it is not hard to see that $W(G,H)\simeq W(K,H\cap K)$.
In particular if $K$ is connected this is the Weyl group of the root
system of $H\cap K$ in $K$.
It is interesting to consider this in the case of real forms of $E_6$
of unequal rank.

\begin{exampleplain}
First  consider $E_6(F_4)$:

\begin{verbatim}
empty: realweyl
Lie type: E6 sc s
(weak) real forms are:
0: e6(f4)
1: e6(R)
enter your choice: 0
there is a unique conjugacy class of Cartan subgroups
Name an output file (return for stdout, ? to abandon):

real weyl group is W^C.((A.W_ic) x W^R), where:
W^C is isomorphic to a Weyl group of type A2
A is trivial
W_ic is a Weyl group of type D4
W^R is trivial
\end{verbatim}
Thus $W(K,H)\simeq W(A_2)\ltimes W(D_4)\simeq S_3\ltimes
W(D_4)$. The action of $S_3$ on $W(D_4)$ is induced by the action of $S_3$
on the Dynkin diagram of $D_4$. By the remarks above we have obtained
the classical
isomorphism 
\begin{equation}
W(F_4)\simeq S_3\ltimes W(D_4).
\end{equation}
\end{exampleplain}

\begin{exampleplain}
Something similar happens for the split real form of $E_6$, in which $K$ 
is of type $C_4$. In this case {\tt realweyl} gives:

\begin{verbatim}
real weyl group is W^C.((A.W_ic) x W^R), where:
W^C is isomorphic to a Weyl group of type A2
A is an elementary abelian 2-group of rank 2
W_ic is a Weyl group of type A1.A1.A1.A1
W^R is trivial
\end{verbatim}

Thus
\begin{equation}
\begin{aligned}
W(C_4)&\simeq S_3\ltimes (\Ztoo^2\times \Ztoo^4)\\
&\simeq S_4\ltimes \Ztoo^4.
\end{aligned}
\end{equation}
\end{exampleplain}

\begin{exampleplain}
\label{ex:complexweyl}
Complex groups.

Fix basic data $(G=G_1\times G_1,\gamma)$  for some group $G_1$,
with $\gamma$ the inner class of the complex group ($\gamma$ switches
the two factors). Thus
$G(\R)=G_1(\C)$ and $K=G_1^\Delta$ (the diagonal embedding).
If $H_1$ is a Cartan subgroup of $G_1$ then
$H=H_1\times H_1$ is a Cartan subgroup of $G$, 
$W(G,H)\simeq W(G_1,H_1)\times
W(G_1,H_1)$, and $W(K,H)=W(G_1,H_1)^\Delta$.

In the notation of the output of {\tt realweyl} all terms are trivial
except {\tt W\^{}C}, which equals $W(G_1,H_1)$.

\begin{verbatim}
empty: type
Lie type: A3.A3 sc C
main: realweyl
there is a unique real form: sl(4,C)
there is a unique conjugacy class of Cartan subgroups
Name an output file (hit return for stdout):

real weyl group is W^C.((A.W_ic) x W^R), where:
W^C is isomorphic to a Weyl group of type A3
A is trivial
W_ic is trivial
W^R is trivial

generators for W^C:
2,5
3,6
1,4
\end{verbatim}
Simple roots $1,2,3$ give the first copy of $G_1$, and $4,5,6$ give
the second. The elements $s_1s_4, s_2s_5$ and $s_3s_6$ generate
$W(G_1,H_1)^\Delta=W(A_3)=S_4$.

\end{exampleplain}

\sec{$K$ orbits on the flag variety and $\Xr$}
\label{s:kgb}

Fix basic data $(G,\gamma)$ (Section \ref{s:basic}) and let $\X^r$ be the reduced parameter space
(Section \ref{s:X}).
Given $x\in\X^r$ let 
$$
\X^r[x]=\{x'\in\X^r\,|\, x'\text{ is $G$-conjugate to }x\}
$$
(see \cite[9.7]{algorithms}).
Choose a preimage $\xi$ of $x$ in $\tX^r$ and let $\theta=\int(\xi)$,
$K=G^\theta$. Then $\X^r[x]$ is isomorphic to $K\bs G/B$. 
In fact, one of the important properties of  $\X^r$ is that it 
captures information about the $K$ orbits on $G/B$
for all $K$ in this inner class. Recall \eqref{e:I} $\{\xi_i\,|\,1\le
i\le n\}$ is a set of
representatives of the strong real forms in this inner class, and for
each $i\in I$ we let $K_i=K_{\xi_i}$.
Then \cite[Corollary 9.9]{algorithms}:

$$
\label{e:allkgb}
\Xr=
\coprod_{i=1}^n\X[x_i]\simeq
\coprod_{i=1}^nK_i\bs G/B.
$$

The structure of $K_i\bs G/B$ is determined to a large extent by the
cross action (cf.~Section \ref{s:X}), and 
{\it Cayley transforms}, which appear in the space $\Xr$ as follows.
Fix $x\in \Xr$ and let $\tau=p(x)\in \I_W$.
Then $\tau$ acts on the roots, and a root $\alpha$ is 
classified as {\it imaginary, real} or {\it complex}
with respect to $\tau$ if $\tau(\alpha)=\alpha,-\alpha$ or neither,
respectively. 
Let  $\xi$ be a preimage of $x$ in $\tXr$. 
If $\alpha$ is imaginary we say it is {\it compact} or {\it
noncompact} with respect to $x$ if this holds with respect to
$\theta_\xi$  (independent of the choice of
$\xi$).
See \cite[Section 14]{algorithms} for details.

Now suppose $\alpha$ is a noncompact imaginary root with respect to $x$. 
Then associated to $\alpha$ and $x$ is a new element of
$\X^r$ denoted $c^\alpha(x)$. It satisfies $p(c^\alpha(x))=s_\alpha
p(x)$, 
and $\alpha$ is real with respect to $s_\alpha\tau$.
Its inverse is the single or
double-valued {\it real} Cayley transform $c_\alpha(x)$:
if $\alpha$ is real with
respect to $x$ then $c_\alpha(x)$ is a set with one or two elements.  

\medskip

The (finite) set of orbits of $K$ on $G/B$ is 
described by the {\tt kgb}
command. 

\begin{exampleplain}
We first consider $SL(2,\R)$ and $PGL(2,\R)$. See \cite[Example 12.20]{algorithms}.

\begin{verbatim}
main: kgb
(weak) real forms are:
0: su(2)
1: sl(2,R)
enter your choice: 1
kgbsize: 3
Name an output file (return for stdout, ? to abandon):
0:  0  0  [n]   1    2
1:  0  0  [n]   0    2
2:  1  1  [r]   2    *  1
\end{verbatim}

Here $G(\R)=SL(2,\R)$, $K(\R)=S^1$ and $K=K(\C)=\C^\times$. 
This acts on $G/B=P^1(\C)=\C\cup\infty$ by $\C^\times\ni
z:w\rightarrow z^2w$.
The entry {\tt kgbsize} gives the number of orbits, $3$ in this case.
These are
$0,\infty$ and $\C^\times$, labelled {\tt \#0,1,2}
respectively
(the first entry in each line) in the output of {\tt kgb}.
(We describe other parts of the output below).
  
Here is the result for $PGL(2,\R)$; see 
\cite[Example 12.25]{algorithms} for more detail.

\begin{verbatim}
0:  0  0  [n]   0    1
1:  1  1  [r]   1    *  1
\end{verbatim}
In this case $G/B$ is the same, but $K=O(2,\C)$; an element from the
non-identity component identifies $0$ and $\infty$, so there are only
two orbits.
\end{exampleplain}

We now explain the other information given in the output. 
Fix a row {\tt i} with corresponding $x_i\in\X^r[x]\simeq K_\xi\bs
G/B$, and write $\mathcal O_i$ for the corresponding orbit of  $K_\xi$
on $G/B$.

The image of $x_i$ in $\I_W$ corresponds to a Cartan subgroup; 
the first entry gives
the length of this
orbit, which is $\dim(\mathcal O_i)-\dim(\mathcal O_0)$. 
In the equal rank case ($\gamma=1$) $\dim(\mathcal O_0)=0$.
The second entry is the number of this Cartan subgroup in the output
of the {\tt cartan} command. 

The simple roots are labelled $1,\dots, n$. 
The term in brackets in row i of the output 
gives the {\it type} of each simple root: {\tt r,C,n,c} for real,
complex, noncompact imaginary, or compact imaginary, respectively.

The next $n$ columns give the cross action
of the simple roots. An entry {\tt j} in column {\tt k} of row {\tt i} means that the
the cross action of the $k^{th}$ simple root takes $x_i$ to $x_j$.

The next $n$ columns give Cayley transforms by noncompact imaginary
roots. There is an entry
in column {\tt k} of row {\tt i} only if the $k^{th}$ simple root is
noncompact imaginary for $x_i$; in this case an entry {\tt j} means
this Cayley transform takes $x_i$ to $x_j$.
These Cayley transforms
are single valued (the inverse, multivalued Cayley transforms are not
listed). 

The final entry is the twisted involution in $W$ corresponding to
$x_i$, 
as a product of
simple reflections (cf.~Section \ref{s:X}).

Here is the $SL(2,\R)$ example again:

\begin{verbatim}
0:  0  0  [n]   1    2
1:  0  0  [n]   0    2
2:  1  1  [r]   2    *  1
\end{verbatim}

Here is the information we can read off from this output.
Orbits {\tt \#0} and {\tt \#1} correspond to the compact Cartan subgroup ({\tt
  \#0}), and are $0$-dimensional. Orbit {\tt \#2}  corresponds to the
split Cartan subgroup and 
is one dimensional. The first two orbits are
interchanged by the cross action of the simple root. For orbits
{\tt \#0,\#1} this root is noncompact, and the Cayley transform takes
each of these orbits to orbit {\tt \#2}.

See \cite[Example 14.19]{algorithms} for a detailed discussion of $Sp(4,\R)$.
\begin{exampleplain}
  
Here is $SL(3,\R)$:

\begin{verbatim}
empty: kgb
Lie type: A2 sc s
there is a unique real form: sl(3,R)
kgbsize: 4
Name an output file (return for stdout, ? to abandon):
0:  0  0  [C,C]   2  1    *  *
1:  1  0  [n,C]   1  0    3  *  2,1
2:  1  0  [C,n]   0  2    *  3  1,2
3:  2  1  [r,r]   3  3    *  *  1,2,1
\end{verbatim}
This is not an equal rank case; the dimension of the closed orbit is
$1$, and the dimension of the unique open orbit is $1+2=3$, the
dimension of $G/B$. 

\end{exampleplain}

\begin{exampleplain}
\label{ex:su22_2}
We illustrate one way to compute the order of $\X^r$ using this
information. Consider equal rank real forms of $SL(4,\C)$, 
i.e. $SU(4), SU(3,1)$ and $SU(2,2)$.

Using {\tt kgb} we compute the number of orbits in each case, given by
{\tt kgbsize}. The result is:

\bigskip

\begin{tabular}{|c|c|}
\hline
$G$&$|K\bs G/B|$\\
\hline
$SU(2,2)$ & 21\\  
\hline
$SU(3,1)$ & 10\\  
\hline
$SU(4,0)$ & 1\\  
\hline
\end{tabular}

\medskip
We will see later (Example \ref{ex:su22_1}) there are two strong real forms
mapping to $SU(4)$ and $SU(3,1)$, and one mapping to $SU(2,2)$. This
means that the set $K\bs G/B$ for $SU(3,1)$, of order $10$, appears twice in
$\Xr$, and similary for $SU(4)$. Therefore the order of $\Xr$ is
\begin{equation}
\label{e:43_2}
2\times 1+2\times 10+1\times 21=43.
\end{equation}
This agrees with the counting done a different way in  
Example \ref{ex:su22_1}.
\end{exampleplain}

\begin{exampleplain}
\label{ex:PSU4_2}
We do the same example again for equal rank forms of the adjoint group
$PSL(4,\C)$.

\medskip

\begin{tabular}{|c|c|}
\hline
$G$&$|K\bs G/B|$\\
\hline
$PU(2,2)$ & 12\\  
\hline
$PU(3,1)$ & 10\\  
\hline
$PU(4,0)$ & 1\\  
\hline
\end{tabular}

\medskip

Since $G$ is adjoint  the map from strong real forms to real forms is
bijective, so $\Xr$ has order $12+10+1=23$.
See Example \ref{ex:PSU4_1}
\end{exampleplain} 

\begin{exampleplain}
\label{ex:kgbsp4}
If $G=Sp(4,\R)$ there are $11$ orbits of $K$ on $G/B$, $4$ of which
are closed.
\begin{verbatim}
 0:  0  0  [n,n]    1   2     6   4
 1:  0  0  [n,n]    0   3     6   5
 2:  0  0  [c,n]    2   0     *   4
 3:  0  0  [c,n]    3   1     *   5
 4:  1  2  [C,r]    8   4     *   *  2
 5:  1  2  [C,r]    9   5     *   *  2
 6:  1  1  [r,C]    6   7     *   *  1
 7:  2  1  [n,C]    7   6    10   *  2,1,2
 8:  2  2  [C,n]    4   9     *  10  1,2,1
 9:  2  2  [C,n]    5   8     *  10  1,2,1
10:  3  3  [r,r]   10  10     *   *  1,2,1,2
\end{verbatim}
For $G=PSp(4,\R)=SO(3,2)$ there are $7$ orbits, $2$ of them closed:
\begin{verbatim}
0:  0  0  [n,n]   0  1    3  2
1:  0  0  [c,n]   1  0    *  2
2:  1  2  [C,r]   5  2    *  *  2
3:  1  1  [r,C]   3  4    *  *  1
4:  2  1  [n,C]   4  3    6  *  2,1,2
5:  2  2  [C,n]   2  5    *  6  1,2,1
6:  3  3  [r,r]   6  6    *  *  1,2,1,2
\end{verbatim}
See Example \ref{ex:sp4block}.
\end{exampleplain}

\begin{exampleplain} Complex Groups.

As in Example \ref{ex:complexweyl} let $G=G_1\times G_1$, 
$K=G_1^\Delta$, 
and $G(\R)=G_1(\C)$. Then 
$B=B_1\times B_1$.
It is well known, and easy to see, that $(g,h)\rightarrow gh\inv$
induces a bijection
\begin{equation}
G_1^\Delta\bs G_1\times G_1/B_1\times B_1
\simeq 
B_1\bs G_1/B_1.
\end{equation}
By  the Bruhat decomposition the right hand side is parametrized by
the Weyl group of $G_1$.

For example take $G(\R)=G_1(\C)=SL(3,\C)$.

\begin{verbatim}
main: type
Lie type: A2.A2 sc C
main: kgb
there is a unique real form: sl(3,C)
kgbsize: 6
Name an output file (hit return for stdout):
0:  0  0  [C,C,C,C]   2  1  2  1    *  *  *  *
1:  1  0  [C,C,C,C]   4  0  3  0    *  *  *  *  2,4
2:  1  0  [C,C,C,C]   0  3  0  4    *  *  *  *  1,3
3:  2  0  [C,C,C,C]   5  2  1  5    *  *  *  *  2,1,3,4
4:  2  0  [C,C,C,C]   1  5  5  2    *  *  *  *  1,2,4,3
5:  3  0  [C,C,C,C]   3  4  4  3    *  *  *  *  1,2,1,3,4,3
\end{verbatim}
Since all roots are complex there are no Cayley transforms. The cross
action of the Weyl group $W(A_2)\simeq S^3$ is simply transitive. 
The last column gives the Weyl group element $(w,w)$ as a product of
simple reflections; just taking the first half of each entry we see
$W=\{id,s_2,s_1,s_2s_1,s_1s_2,s_1s_2s_1\}$.
See Example \ref{ex:complexweyl}.

\end{exampleplain}

\sec{Representation Theory: Blocks and the two-sided parameter space}
\label{s:rep}
The final ingredient in the algorithm is the {\it two-sided parameter
space}. 
Let $\chG$ be the dual group of $G$. The involution $\gamma\in
\Out(G)$ defines an involution $\ch\gamma\in\Out(\chG)$
 and $(\chG,\ch\gamma)$ is
also  basic data.
See \cite[Section 2]{algorithms} for details.
Let $\ch\X$ be the one-sided parameter space
defined by $(\chG,\ch\gamma)$. 

Suppose $x\in \X$ and  $\xi\in\tX$ is a pre-image of $x$ in
$\tX$. Then $\theta_\xi$ restricted to $H$ is independent of the
choice of $\xi$, and we denote it $\theta_{x,H}$. We also write
$\theta_{x,H}\in \End(\h)$ for its differential. There is a natural
pairing between $\h$ and $\ch\h$ (a Cartan subalgebra on the dual
side); the adjoint $\theta_{x,H}^t$ of $\theta_{x,H}$ is an element of $\End(\ch\h)$. 

We can now define the {\it two-sided parameter space}:
\begin{equation}
\label{e:caZ}
\caZ=\{(x,y)\in\X\times\ch\X\,|\, \theta_{x,H}^t=-\theta_{y,\ch H}\}.
\end{equation}
See \cite[Section 10]{algorithms}.
We define the {\it reduced} two-sided parameter space $\caZr$ by replacing
$\X$ and $\ch\X$ with the corresponding reduced spaces (cf.~Section \ref{s:X}).
This is a finite set.

With the obvious notation $W\times \ch W$ acts on $\caZ$ and $\caZr$;
this is the
cross action. By the condition relating $x$ and $y$ in \eqref{e:caZ},
if $\alpha$ is a real root with 
respect to  $x$ then
$\ch\alpha$ is an imaginary root with respect to  $y$, and vice-versa.
Thus Cayley
transforms are defined on $\caZr$, via a real root on one side, and an
imaginary noncompact root on the other.

The main result in \cite{algorithms}, Theorem 10.3, 
(also see \cite[Theorem 7.15]{algorithms}),
is that the space $\caZ$
parametrizes the irreducible representations of strong real forms of
$G$ ``up to translation'';
equivalently  with certain regular integral infinitesimal characters. 
Because of our emphasis on the reduced parameter
space, we give a slightly different version here. 

Apply the construction in the paragraph preceding \eqref{e:X} to
choose
a subset $Z^{\vee r}$ of the center $\ch Z$ of $\ch G$.
The pairing of $\h$ and $\ch\h$ gives an isomorphism $\h^*\simeq \ch\h$, and let
\begin{equation}
L=\{\lambda\in \h^*\,|\,\exp(2\pi i\lambda)\in Z^{\vee r}\}.
\end{equation}
Note that $L$ is a subset of the weight lattice
$P=\{\lambda\in\h^*\,|\, \langle\lambda,\ch\alpha\rangle\in\Z\}$. 
Choose a set $\Lambda$ of representatives of $L/X^*(H)$,
satisfying
$\langle \lambda,\ch\alpha\rangle\ne 0$ for all
$\lambda\in\Lambda$ and all roots $\alpha$. This is a finite set.

Choose a set $\{\xi_i\,|\, 1\le i\le n\}$  as in \eqref{e:I} 
For each $i$ let $G_i(\R)$ be the real form of $G$
defined by $\xi_i$, and let $\Pi(G_i(\R),\Lambda)$ be the set of
irreducible admissible  representations of $G_i(\R)$ with
infinitesimal character contained in $\Lambda$.

\begin{theorem}
\label{t:main}
There is a natural bijection

\begin{equation}
\caZ^r\overset{1-1}\longleftrightarrow \coprod_{i=1}^n\Pi(G_i(\R),\Lambda).
\end{equation}
\end{theorem}

This differs from \cite[Theorem 10.3]{algorithms} in several ways. First
of all on the left hand side we're using $\caZ^r$ instead of $\caZ$.
On the right hand side the union is over
strong real forms in the sense of the reduced parameter space,
instead
of the possibly infinite set $\I/G$ of \cite[5.15]{algorithms}.
Finally 
the set $\Lambda=L/X^*(H)$  is contained in the
corresponding set $P/X^*(H)$ of \cite{algorithms}.

\medskip

Fix $(x,y)\in \caZr$, and consider the representations associated to
the pairs $(x',y')$ in the subset $\X^r[x]\times \X^{\vee r}[y]$ of $\caZr$. Choose $i$ so that
$x$ is $G$-conjugate to $x_i$. Then these are representations of
$G_i(\R)$, all with the same infinitesimal character. In fact this set
of irreducible representations is a {\it block} in the sense of
\cite[Chapter 9]{green}, and every block is obtained this way.

The entire construction is obviously symmetric in $G$ and $\ch G$, so
an element $(x,y)\in\caZ^r$, which defines a representation $\pi$ of a
real form of $G$, also defines a representation $\ch\pi$ of a
real form of $\ch G$. The  map $\pi\rightarrow \ch\pi$ is a version of
{\it Vogan Duality}; see \cite{ic4} and 
\cite[1.35 and Corollary 10.9]{algorithms}. In particular this is a duality of
{\it blocks}: the set $\X^r[x]\times\X^{\vee r}[y]$ defines a block $\caB$
of a real form of $G$, a block $\ch\caB$ of a real form of $\ch G$,
and gives a bijection $\caB\leftrightarrow\ch\caB$.



\medskip

Before we look at individual representations, it is helpful to look at
the sizes of blocks, which are given by the {\tt blocksizes} command.

\begin{exampleplain}
\label{ex:A1blocksizes}
Here are the blocks for real forms of $SL(2,\C)$:
\begin{verbatim}
main: type
Lie type: A1 sc s
main: blocksizes
  0  1
  1  3
\end{verbatim}
The rows and columns correspond to real forms of $G$ and $\ch G$, respectively.
Adding some labelling by hand the picture is:

\begin{verbatim}
        SO(3) SO(2,1)
SU(2)   0     1
SL(2,R) 1     3
\end{verbatim}
Thus $SU(2)$ has a single block, dual to
a block of $SO(2,1)$, and
$SL(2,\R)$ has two blocks, dual to $SO(3)$, and $SO(2,1)$,
respectively.
See Example \ref{ex:sl2block}.
\end{exampleplain}

\begin{exampleplain}
\label{ex:psl2bocksizes}
The corresponding output for $PSL(2,\C)$ is the same, but there is a
subtle point here.

\begin{verbatim}
main: type
Lie type: A1 sc s
main: blocksizes
  0     1
  1     3
\end{verbatim}
Since the dual group $SL(2,\C)$ is not adjoint, the map from strong
real forms to real forms (on the dual side) is not injective. There
are two strong real forms mapping to $SU(2)$, and we can label the
strong real forms $SU(2,0), SU(1,1)$ and $SU(0,2)$. 
Consequently there are
two blocks of $SO(2,1)$ of size $1$. 
See Example
 \ref{ex:blockso21}.
We could display this information by constructing a table by hand
showing {\it strong} real forms on the dual side:
\begin{verbatim}
        SU(2,0) SU(0,2) SL(2,R)
SO(3)   0       0       1
SO(2,1) 1       1       3
\end{verbatim}
\end{exampleplain}

Blocks
 tend to be concentrated on the quasisplit
forms of $G$ or $\ch G$.

\begin{exampleplain}
Here is the output of {\tt blocksizes} for real forms of  $Sp(12,\C)$, with the
real forms added:

\begin{verbatim}
         SO(13) SO(12,1) SO(11,2) SO(10,3) SO(9,4) SO(8,5) SO(7,6)
Sp(6)    0      0        0        0        0       0       1
Sp(5,1)  0      0        0        0        0       0       36
Sp(4,2)  0      0        0        0        0       0       315
Sp(3,3)  0      0        0        0        0       0       680
Sp(12,R) 1      13       108      556      1975    4707    7416
\end{verbatim}
\end{exampleplain}

\begin{exampleplain}
\label{ex:so12}
In the case of equal rank  real forms of $SO(12,\C)$ 
we get some other elements:

\begin{verbatim}
         SO(12) SO(10,2) SO*(12) SO*(12) SO(8,4) SO(6,6)
SO(12)   0      0        0       0       0       1
SO(10,2) 0      0        0       0       15      66
SO*(12)  0      0        0       60      0       692
SO*(12)  0      0        60      0       0       692
SO(8,4)  0      15       0       0       300     885
SO(6,6)  1      66       692     692     885     2320
\end{verbatim}
Note that this diagram is symmetric: $SO(12,\C)$ is self-dual (and the 
equal rank and split inner classes coincide).
For a discussion of the two versions of $SO^*(12)$ see
Example \ref{ex:SO12_2}. Note that the two versions of $SO^*(12)$ on the dual side
account for the fact that $SO(6,6)$ has two distinct
blocks of size $692$.
\end{exampleplain}

\begin{exampleplain}
We break the symmetry between $G$ and $\ch G$ 
of the previous example by taking $Spin(12,\C)$,
which is dual to $PSO(12,\C)$.

\begin{verbatim}
           PSO(12) PSO(10,2) PSO*(12) PSO*(12) PSO(8,4) PSO(6,6)
Spin(12)   0       0         0        0        0        1
Spin(10,2) 0       0         0        0        15       87
Spin*(12)  0       0         0        60       0        692
Spin*(12)  0       0         60       0        0        692
Spin(8,4)  0       15        0        0        300      915
Spin(6,6)  1       66        436      436      885      2180
\end{verbatim}

There is an interesting phenomenon here. Note that $Spin(6,6)$ has a
block $\caB_{Spin}$ of size $66$, dual to a block of
$PSO(10,2)$. From the previous example
$SO(6,6)$ also has a block of size $66$, denoted $\caB_{SO}$, 
dual to a block for
$SO(10,2)$. From the output of the {\tt block} command one can see 
these blocks are isomorphic.

All representations in a block have the same central
character, and the representations in $\caB_{Spin}$ factor 
to the image of $Spin(6,6)$ in $SO(6,6)$. 
The image of this map
has index $2$ in $SO(6,6)$. Therefore it is not obvious that
$\caB_{Spin}$ and $\caB_{SO}$ should be isomorphic (and this doesn't
happen for the blocks of $Spin(6,6)$ of size $436$ and $2,180$).

The explanation is seen by looking at the dual side. The map
$SO(10,2)\rightarrow PSO(10,2)$ {\it is} surjective. It follows 
that these dual blocks are isomorphic, hence the blocks themselves are
isomorphic. 
\end{exampleplain}

We leave it to the reader to see that the symmetry between the two
versions of $SO^*(12)$ is broken by taking the ``non-diagonal''
quotient of $SO(12,\C)$.
See Examples 
\ref{ex:SO8nondiag} and
\ref{ex:SO12_2}.

\begin{exampleplain}
Here is the output of {\tt blocksizes} for $E_8$, with the real forms
added:

\begin{verbatim}
              compact quaternionic split
compact       0       0            1
quaternionic  0       3150         73410
split         1       73410        453060
\end{verbatim}
The quaternionic real form is the one denoted {\tt e8(e7.su(2))}
in the software.
The software computes the structure of the block of size 
$453,060$ easily. However computing Kazhdan-Lusztig-Vogan polynomials
for this block required special techniques. See {\tt
www.liegroups.org} for more information.
\end{exampleplain}
\bigskip

We now look at individual representations in a few examples.
The cases of real forms of $SL(2,\C)$ and $PSL(2,\C)$ are explained 
in detail in \cite[Example 12.20]{algorithms}. 
In particular see the table at the end of Section 12.
Here is a brief
summary. 

\begin{exampleplain}
\label{ex:sl2block}
$SL(2,\R)$ and $SU(2)$:
\begin{verbatim}
empty: block
Lie type: A1 sc s
(weak) real forms are:
0: su(2)
1: sl(2,R)
enter your choice: 1
possible (weak) dual real forms are:
0: su(2)
1: sl(2,R)
enter your choice: 1
Name an output file (return for stdout, ? to abandon):
0(0,1):  0  0  [i1]  1   (2,*)
1(1,1):  0  0  [i1]  0   (2,*)
2(2,0):  1  1  [r1]  2   (0,1)   1
\end{verbatim}
This is the block of $SL(2,\R)$ of size $3$ of example \ref{ex:A1blocksizes}.
We can take the infinitesimal character to be $\rho$,
and this block contains the two discrete series  representations {\tt
  \#0,\#1}, and
the trivial representation {\tt \#2}. The singleton block of $SL(2,\R)$ is
found taking the dual group to be compact:

\begin{verbatim}
possible (weak) dual real forms are:
0: su(2)
1: sl(2,R)
enter your choice: 0
Name an output file (return for stdout, ? to abandon):
0(2,0):  1  1  [rn]  0   (*,*)   1
\end{verbatim}
This is an irreducible (non-spherical) principal series representation of $SL(2,\R)$
at $\rho$ (with odd $K$-types).
\end{exampleplain}

\begin{exampleplain}
\label{ex:blockso21}
The adjoint group
$PSL(2,\R)\simeq SO(2,1)$ has a block of size $3$ dual to the
preceding one:

\begin{verbatim}
0(0,2):  0  0  [i2]  0   (1,2)
1(1,0):  1  1  [r2]  2   (0,*)   1
2(1,1):  1  1  [r2]  1   (0,*)   1
\end{verbatim}
We can take infinitesimal character $\rho$, and the block consists of 
the unique discrete series representation {\tt \#0} and the two one-dimensional
representations {\tt \#1,\#2}. These are all of the irreducible representations of
$PSL(2,\R)$ with infinitesimal character $\rho$.

As before taking the dual group to be compact we obtain a singleton:

\begin{verbatim}
possible (weak) dual real forms are:
0: su(2)
1: sl(2,R)
enter your choice: 0
Name an output file (return for stdout, ? to abandon):
0(1,0):  1  1  [rn]  0   (*,*)   1  
\end{verbatim}

This block occurs at infinitesimal character $2\rho$; in this case the
set $\Lambda$ of Theorem \ref{t:main} can be taken to be
$\{\rho,2\rho\}$.  Note that $SO(2,1)$ has two irreducible
representation $PS_\pm$ at infinitesimal character $2\rho$.  The fact
that there are two such representations (each of which is a block) is
reflected in the fact that on the dual side there are two strong real
forms $SU(2,0$ and $SU(0,2)$. See Examples \ref{ex:A1blocksizes} and
\ref{ex:psl2bocksizes}, and the table at the end of \cite[Section
12]{algorithms}.
\end{exampleplain}

We illustrate the other information in the output of {\tt block} by
looking at $Sp(4,\R)$.
See \cite[Example 14.19]{algorithms}.

\begin{exampleplain}
\label{ex:sp4block}
Here is the block of size $12$ for $Sp(4,\R)$:  
\begin{verbatim}
 0( 0,6):  0  0  [i1,i1]   1   2   ( 6, *)  ( 4, *)
 1( 1,6):  0  0  [i1,i1]   0   3   ( 6, *)  ( 5, *)
 2( 2,6):  0  0  [ic,i1]   2   0   ( *, *)  ( 4, *)
 3( 3,6):  0  0  [ic,i1]   3   1   ( *, *)  ( 5, *)
 4( 4,4):  1  2  [C+,r1]   8   4   ( *, *)  ( 0, 2)   2
 5( 5,4):  1  2  [C+,r1]   9   5   ( *, *)  ( 1, 3)   2
 6( 6,5):  1  1  [r1,C+]   6   7   ( 0, 1)  ( *, *)   1
 7( 7,2):  2  1  [i2,C-]   7   6   (10,11)  ( *, *)   2,1,2
 8( 8,3):  2  2  [C-,i1]   4   9   ( *, *)  (10, *)   1,2,1
 9( 9,3):  2  2  [C-,i1]   5   8   ( *, *)  (10, *)   1,2,1
10(10,0):  3  3  [r2,r1]  11  10   ( 7, *)  ( 8, 9)   1,2,1,2
11(10,1):  3  3  [r2,rn]  10  11   ( 7, *)  ( *, *)   1,2,1,2
\end{verbatim}

\medskip

The $12$ representations, numbered $0,\dots,11$, 
are parametrized by pairs $(x,y)$, the second entry on each
line, from the corresponding {\tt kgb} commands for $G$ and $\ch G$.
See Example \ref{ex:kgbsp4}; note that for $G=Sp(4,\R)$ there are $12$ orbits
numbered {\tt 0,\dots, 11} and for $\ch G=SO(3,2)$ there are $7$ orbits,
labelled {\tt 0,\dots,6}.

The next two columns give the Cartan subgroup, and length of the
parameter, just as in the output of {\tt kgb} for $x$ (Section \ref{s:kgb}).

The term in brackets list each of the simple roots as:

\begin{itemize}
\item compact imaginary: {\tt ic}
\item noncompact imaginary, type I/II: {\tt i1/i2}
\item complex: {\tt C+,C-}
\item real, not satisfying the parity condition: {\tt rn}
\item real, satisfying the parity condition type I/II: {\tt r1/r2}
\end{itemize}
For information about type I/II roots and the parity condition see \cite[Section 8.3]{green}.

The next two columns give cross actions of the simple roots, similar 
to the {\tt kgb} command, followed by two columns for 
Cayley transforms. These follow from the cross action/Cayley
transforms for {\tt kgb} on both the $G$ and $\ch G$ side. 
In this setting (unlike {\tt kgb}) these can be double
valued, even for a noncompact imaginary root, since this corresponds 
to a real root on the dual side.

The final column is the corresponding twisted involution, exactly as
in the {\tt kgb} command for $x$.
\end{exampleplain}

\begin{exampleplain}
\label{ex:e6}
We conclude with a large example.

\begin{verbatim}
empty: block
Lie type: E6 sc c
(weak) real forms are:
0: e6
1: e6(so(10).u(1))
2: e6(su(6).su(2))
enter your choice: 2
possible (weak) dual real forms are:
0: e6(f4)
1: e6(R)
enter your choice: 0
\end{verbatim}

This is the block for $E_6(A_5\times A_1)$ (simply connected), dual to
$E_6(F_4)$ (adjoint).

{\scriptsize
\begin{verbatim}
 0( 851,44):   8  4  [C+,rn,rn,rn,rn,C+]   2   0   0   0   0   1
 1(1013,43):   9  4  [C+,rn,rn,rn,C+,C-]   4   1   1   1   3   0
 2(1014,42):   9  4  [C-,rn,C+,rn,rn,C+]   0   2   5   2   2   4
 3(1165,41):  10  4  [C+,rn,rn,C+,C-,rn]   8   3   3   6   1   3
 4(1166,40):  10  4  [C-,rn,C+,rn,C+,C-]   1   4   9   4   8   2
 5(1167,39):  10  4  [rn,rn,C-,C+,rn,C+]   5   5   2   7   5   9
 6(1304,38):  11  4  [C+,C+,C+,C-,rn,rn]  11  10  10   3   6   6
 7(1305,37):  11  4  [rn,C+,rn,C-,C+,C+]   7  14   7   5  14  12
 8(1306,36):  11  4  [C-,rn,C+,C+,C-,rn]   3   8  13  11   4   8
 9(1307,35):  11  4  [rn,rn,C-,C+,C+,C-]   9   9   4  12  13   5
10(1430,34):  12  4  [C+,C-,C-,rn,rn,rn]  18   6   6  10  10  10
11(1431,33):  12  4  [C-,C+,C+,C-,rn,rn]   6  18  15   8  11  11
12(1432,32):  12  4  [rn,C+,rn,C-,C+,C-]  12  19  12   9  17   7
13(1433,31):  12  4  [rn,rn,C-,C+,C-,rn]  13  13   8  16   9  13
14(1434,30):  12  4  [rn,C-,rn,rn,C-,C+]  14   7  14  14   7  19
15(1536,29):  13  4  [C+,C+,C-,C+,rn,rn]  24  24  11  20  15  15
16(1537,28):  13  4  [rn,C+,C+,C-,C+,rn]  16  22  20  13  21  16
17(1538,27):  13  4  [rn,C+,rn,C+,C-,C+]  17  23  17  21  12  23
18(1539,26):  13  4  [C-,C-,C+,rn,rn,rn]  10  11  24  18  18  18
19(1540,25):  13  4  [rn,C-,rn,rn,C+,C-]  19  12  19  19  23  14
20(1620,24):  14  4  [C+,C+,C-,C-,C+,rn]  28  26  16  15  25  20
21(1621,23):  14  4  [rn,C+,C+,C-,C-,C+]  21  29  25  17  16  27
22(1622,22):  14  4  [rn,C-,C+,rn,C+,rn]  22  16  26  22  29  22
23(1623,21):  14  4  [rn,C-,rn,C+,C-,C-]  23  17  23  27  19  17
24(1624,20):  14  4  [C-,C-,C-,C+,rn,rn]  15  15  18  28  24  24
25(1684,19):  15  4  [C+,C+,C-,rn,C-,C+]  32  31  21  25  20  30
26(1685,18):  15  4  [C+,C-,C-,C+,C+,rn]  34  20  22  34  31  26
27(1686,17):  15  4  [rn,C+,C+,C-,rn,C-]  27  33  30  23  27  21
28(1687,16):  15  4  [C-,C+,rn,C-,C+,rn]  20  34  28  24  32  28
29(1688,15):  15  4  [rn,C-,C+,C+,C-,C+]  29  21  31  33  22  33
30(1730,14):  16  4  [C+,C+,C-,rn,rn,C-]  37  36  27  30  30  25
31(1731,13):  16  4  [C+,C-,C-,C+,C-,C+]  38  25  29  35  26  36
32(1732,12):  16  4  [C-,C+,rn,rn,C-,C+]  25  38  32  32  28  37
33(1733,11):  16  4  [rn,C-,C+,C-,rn,C-]  33  27  36  29  33  29
34(1734,10):  16  4  [C-,C-,rn,C-,C+,rn]  26  28  34  26  38  34
35(1760, 9):  17  4  [C+,rn,C+,C-,C+,C+]  40  35  39  31  40  39
36(1761, 8):  17  4  [C+,C-,C-,C+,rn,C-]  41  30  33  39  36  31
37(1762, 7):  17  4  [C-,C+,rn,rn,rn,C-]  30  41  37  37  37  32
38(1763, 6):  17  4  [C-,C-,rn,C+,C-,C+]  31  32  38  40  34  41
39(1777, 5):  18  4  [C+,rn,C-,C-,C+,C-]  43  39  35  36  42  35
40(1778, 4):  18  4  [C-,rn,C+,C-,C-,C+]  35  40  42  38  35  43
41(1779, 3):  18  4  [C-,C-,rn,C+,rn,C-]  36  37  41  43  41  38
42(1786, 2):  19  4  [C+,rn,C-,rn,C-,C+]  44  42  40  42  39  44
43(1787, 1):  19  4  [C-,rn,C+,C-,C+,C-]  39  43  44  41  44  40
44(1790, 0):  20  4  [C-,rn,C-,rn,C-,C-]  42  44  43  44  43  42             
\end{verbatim}}

The last column (twisted involution) has been deleted, as have the
Cayley transforms, which are all empty, i.e. {\tt (*,*)}.

All of these
representations come from the same Cartan subgroup. This is
the entry {\tt 4} in each row;  by the output of the {\tt cartan}
command this is the most split Cartan subgroup, and is isomorphic to
$\R^\times\times\R^\times\times\C^\times\times\C^\times$. 
Note that every root is of type {\tt C$\pm$} or {\tt rn}, so there are
no Cayley transforms. This is surprising, since this 
group has $5$ conjugacy classes of Cartan subgroups
including a compact Cartan subgroup but not a split one.
The reason is that the dual
block is for the adjoint group $E_6(F_4)$, which has only one
conjugacy class of Cartan subgroup (cf. Example \ref{ex:E6f4}).

There are $1,791$ orbits of $K$ on $G/B$, and $45$ for the dual group,
given by the {\tt kgb} command. In the parameters $(x,y)$ given in
parentheses each $0\le y\le 44$ appears once,
since the dual Cartan subgroup is connected. The fact that $45$
distinct values of $x$ (between $0$ to $1790$) occur is more subtle.
\end{exampleplain}

\sec{The Geometry of $\Xr$}
\label{s:geometryofX}

Fix $(G,\gamma)$ and let $\Xr$ be the reduced parameter space defined in
Section \ref{s:X}.
It is helpful to give some detail about the structure of $\Xr$.
Recall (Section~\ref{s:X}) there is a map from $\X^r$ to the space $\I_W$ of twisted
involutions in $W$. We discuss the fibers of
this map.

\subsec{The Adjoint Case}
\label{s:adjoint}

Fix $(G,\gamma)$ with $G$ adjoint.
This case is simplest,
since real forms and strong real forms
coincide.  Fix $\tau\in \I_W$.

Let $H^{-\tau}=\{h\in H\,|\, \tau(h)=h\inv\}$. 
There is a natural simply transitive action of
$H^{-\tau}/(H^{-\tau})^0$ on $\X^r_\tau$ \cite[Section 11]{algorithms}.
Note that 
$H^{-\tau}/(H^{-\tau})^0\simeq\Ztoo^b$ where $b$ is the number of
$S^1$ factors of the real form of $H$ defined by $\tau$.
Fix $x\in \X^r_\tau$.
Via a choice of basepoint in $\X^r_\tau$ the output of {\tt cartan} identifies $\X^r_\tau$ with $\{0,1,\dots,
2^b-1\}$.  The element of $\Ztoo^b$ corresponding to $0\le k\le 2^b-1$ is
the binary expansion of $k$.

The rank of $\X^r_\tau$ is also given by {\tt fiber rank} in the output
of {\tt cartan}. For each Cartan the number {\tt \#X\_r} is the number
of strong involutions lying over this conjugacy class in $\I_W$, i.e. 
{\tt twisted involution orbit size}$\times$2\^{\tt fiber rank}.

The imaginary Weyl group $W_{i,\tau}$
(the Weyl group of the imaginary roots with
respect to $\tau$)
 acts on $\X^r_\tau$, and the orbits are
in one to one correspondence with real forms  containing this Cartan
subgroup. The {\tt cartan} command gives this decomposition.
For example a line 

\smallskip
\noindent{\tt real form \#2: [0,1,2] (3)} 
\smallskip

\noindent means that corresponding to
real form {\tt \#2} (from the output of the {\tt realform} command) are
three elements of $\X^r_\tau$, labelled {\tt 0,1,2} (the number in
parentheses is the number of elements in brackets). 

\begin{exampleplain}
\label{ex:PSU4_1}
Here are the equal rank real forms of $PSL(4,\C)$:

\begin{verbatim}
empty: cartan
Lie type: A3 ad c
(weak) real forms are:
0: su(4)
1: su(3,1)
2: su(2,2)
enter your choice: 2
Name an output file (return for stdout, ? to abandon):

Cartan #0:
split: 0; compact: 3; complex: 0
canonical twisted involution:
twisted involution orbit size: 1;  fiber rank: 3;  #X_r: 8
imaginary root system: A3
real root system is empty
complex factor is empty
real form #2: [0,2,5] (3)
real form #1: [1,3,4,6] (4)
real form #0: [7] (1)
\end{verbatim}

This is the fundamental fiber $\tau=\delta$. In this case $W_{i,\delta}=W$. 
This fiber has $2^3=8$ elements. There is a fixed point of the action
of $W$ on $\X_\delta$, corresponding to the compact real form $PSU(4)$
({\tt realform \#0}).
The largest orbit, of size $4$, corresponds to the group $SU(3,1)$. 

Here are the remaining Cartan subgroups:

\begin{verbatim}
Cartan #1:
split: 0; compact: 1; complex: 1
canonical twisted involution: 1,2,3,2,1
twisted involution orbit size: 6;  fiber rank: 1;  #X_r: 12
imaginary root system: A1
real root system: A1
complex factor is empty
real form #2: [0] (1)
real form #1: [1] (1)


Cartan #2:
split: 1; compact: 0; complex: 1
canonical twisted involution: 2,1,3,2
twisted involution orbit size: 3;  fiber rank: 0;  #X_r: 3
imaginary root system is empty
real root system: A1.A1
complex factor: A1
real form #2: [0] (1)
\end{verbatim}

There are $6$ elements $\tau$ here
for which $\X^r_\tau$ has order $2$; both noncompact
real forms appear. 
Finally only the quasisplit group
$PSU(3,3)$ contains the maximally split Cartan, in which case there
are $3$ twisted involutions $\tau$, each with $|\X^r_\tau|=1$.

Therefore $\X^r$ has $1\times 8+6\times 2+3\times 1=23$ elements.
See Example \ref{ex:PSU4_2}
\end{exampleplain}

\subsec{Strong Real Forms and Fibers}
\label{s:fiber}

The {\tt strongreal} command gives information about strong real
forms (see Section \ref{s:X}), which may also be interpreted as information about the fibers
$p:\Xr\rightarrow \I_W$ in the general case.


The strong real forms ``containing'' $(H,\tau)$ are parametrized by
$\X^r_\tau/W_{i,\tau}$.
In particular every strong real form contains the fundamental Cartan
subgroup $(H,\delta)$, so strong real forms are parametrized by 
$\X^r_{\delta}/W_{i,\delta}$.
Things are particularly simple if $G$ has equal rank, say $n$, in which case
\begin{equation}
\X^r_\delta/W\simeq \{h\in H\,|\, h^2\in Z^r\}/W.
\end{equation}
Note that $|\X^r_\delta|=|Z^r|2^n$ in this case.

In general not every element of $Z^r$ is of the form $x^2$ for
some $x\in\X$. 
The software labels
the elements  which are of this form  as $\{z_0,\dots, z_{r-1}\}$, denoted {\tt class \#0, class
\#1,\dots class \#r-1}. In fact $r=2^s$ for some $s$.
For each $z_j$ the set 
$\X^r_\tau(z_j)=\{x\in \X_\tau\,|\, x^2=z_j\}$
is isomorphic to $H^{-\tau}/(H^{-\tau})^0$ (independent of $j$) which
has cardinality $2^b$
($b$ is the number of $S^1$ factors as in the preceding section)
so
$|\Xr_\tau|=2^{b+s}$.

The  group $W_{i,\tau}$ acts on $\X_\tau(z_j)$, and the orbits
correspond to strong real forms.
All of this information is given by the {\tt strongreal} command. 
See \cite[Section 11]{algorithms} and the {\tt help} file for 
the {\tt strongreal} command.

Note that the order of $\X^r$ may be computed by summing over
Cartan subgroups, taking into account the size of the conjugacy
classes of twisted involutions.

\begin{exampleplain}
\label{ex:su22_1}
Consider the equal rank inner class of $SL(4,\C)$.
Compare  Example \ref{ex:PSU4_1}.

\begin{verbatim}
main: strongreal
Lie type: A3 sc c
(weak) real forms are:
0: su(4)
1: su(3,1)
2: su(2,2)
enter your choice: 0
there is a unique conjugacy class of Cartan subgroups
Name an output file (hit return for stdout):

there are 2 real form classes:

class #0:
real form #2: [0,1,2,4,5,6] (6)
real form #0: [3] (1)
real form #0: [7] (1)

class #1:
real form #1: [0,2,3,5] (4)
real form #1: [1,4,6,7] (4)
\end{verbatim}

First of all, reading the lines beginning {\tt realform}, we see there
is one strong real form mapping to $SU(2,2)$ (real form {\tt \#2}),
two mapping to $SU(4)$ (real form {\tt \#0}), and two
mapping to $SU(3,1)$ (real form {\tt \#1}). We can think of these as 
$SU(4,0), SU(3,1), SU(2,2), SU(1,3)$, and $SU(0,4)$.

This is the fundamental fiber, i.e. $\tau=\delta$.
In this case $Z\simeq \Z_4$ and $Z/Z^2\simeq \Ztoo$, so $|Z^r|=2$.
In this case $b$ (the number of $S^1$ factors) is $3$, so
$|\X^r_\delta(z)|=8$ for each $z\in Z^r$ (it is never $0$ in this case), 
and $|\X^r_\delta|=16$.
We can take $Z/Z^2=\{z_0,z_1\}=\{I,iI\}$. These are the elements {\tt class
\#0} and {\tt class \#1} respectively.

There are three orbits of $W=W_{i,\delta}$ on $\X^r_\delta(I)$.
The $W$-orbit of $\diag(1,1,-1,-1)$ has $6$ elements; given in the
line 
{\tt real form \#2: [0,1,2,4,5,6] (6)}; the corresponding  real form  is
$SU(2,2)$.
The elements $\pm I$ are each fixed by $W$, and correspond to the two lines
beginning {\tt real form \#0}. (The software does not specify which of these is
$I$, and which $-I$).

We could have chosen $z_0=-I$ instead of $I$.
Then $SU(2,2)$ would be given by $i$ times the preceding
ones, i.e. the $W$-orbit of $\diag(i,i,-i,-i)$. Similarly the two
strong real forms mapping to $SU(2)$ would be $\pm iI$.
The software does not make an actual choice of $Z^r$; the
output of the software, and the combinatorics of the algorithm, are
independent of any such choice. See \cite[Section 13]{algorithms}.

Now consider $z=iI$.
Then $\X^r_\delta(iI)$ consists of 
the elements $\zeta(\pm 1,\dots,
\pm1)$ where $\zeta=e^{\pi i/4}$, with an odd number of plus
signs. The $8$ such elements constitute $2$ $W$-orbits, hence the two
strong real forms mapping to $SU(3,1)$.
\end{exampleplain}

Here are the strong real forms containing the other two Cartan
subgroups in this example.

\begin{verbatim}
real: strongreal
cartan class (one of 0,1,2): 1
Name an output file (return for stdout, ? to abandon):

there are 2 real form classes:

class #0:
real form #2: [0,1] (2)

class #1:
real form #1: [0] (1)
real form #1: [1] (1)
\end{verbatim}
The compact real form {\tt real form \#0} only contains the compact
Cartan subgroup and doesn't occur.
In this case $|\X^r_\tau(I)|=|\X^r_\tau(iI)|=2$, 
$|\X^r_\tau|=4$, and 
both $SU(2,2)$ and $SU(3,1)$ occur. 

\begin{verbatim}
real: strongreal
cartan class (one of 0,1,2): 2
Name an output file (return for stdout, ? to abandon):

real form #2: [0] (1)
\end{verbatim}
For the most split Cartan subgroup  $|\X^r_\tau(I)|=0$, $|\X^r_\tau(iI)|=1$, and 
only $SU(2,2)$ occurs.

The number of twisted involutions in each conjugacy class are $1,6$
and $3$, as given by the {\tt cartan} command. Therefore the
cardinality of $\X^r$ is
\begin{equation}
\label{e:43_1}
16\times 1+4\times 6+1\times 3=43.
\end{equation}
See Example 
\ref{ex:su22_2}.

\begin{exampleplain}
Type $A_{2n}$ is a little different than $A_{2n+1}$. Here are the
fundamental fiber and the equal rank strong real forms of  $SL(5,\C)$:
\begin{verbatim}
Lie type: A4 sc c
main: strongreal
(weak) real forms are:
0: su(5)
1: su(4,1)
2: su(3,2)
enter your choice: 0
there is a unique conjugacy class of Cartan subgroups
Name an output file (hit return for stdout):

real form #2: [0,1,2,4,5,6,8,10,11,13] (10)
real form #1: [3,9,12,14,15] (5)
real form #0: [7] (1)
\end{verbatim}
In this case $Z\simeq \Z/5\Z$ and $Z/Z^2=1$. We can take $z=I$ in this
case, and the map from strong real forms to real forms is bijective.
\end{exampleplain}

We conclude with a discussion of 
the inner class of $SL(n,\R)$.
If $n$ is odd this is the
unique real form in this inner class; if $n$ is even there is one
other real form $SL(n/2,\mathbb H)$. 
We can take $\theta(\diag(z_1,\dots, z_n))=\diag(\frac1{z_n},\dots,
\frac1{z_1})$, and 
$Z^r=\{I\}$ if $n$ is odd, or $\{\pm I\}$ if $n$ is even.

Here is an example of each case.

\begin{exampleplain}
There  is only one strong real form in
the inner class of $SL(3,\R)$:
\begin{verbatim}
empty: strongreal
Lie type: A2 sc s
there is a unique real form: sl(3,R)
\end{verbatim}
Here $Z^r=\{I\}$ and 
\begin{equation}
\X_{\delta}(I)=\{\diag(z,w,z)\,|\, z^2w=1\}/
\{\diag(ac,b^2,ac)\,|\,abc=1\}=I
\end{equation}
\end{exampleplain}

\begin{exampleplain}
Now consider the inner class of $SL(4,\R)$, which also contains the
real form $SL(2,\mathbb H)$.
In this case $Z=\pm I$.

\begin{verbatim}
empty: strongreal
Lie type: A5 sc s
(weak) real forms are:
0: sl(3,H)
1: sl(6,R)
enter your choice: 0
there is a unique conjugacy class of Cartan subgroups
Name an output file (return for stdout, ? to abandon):

there are 2 real form classes:

class #0:
real form #1: [0,1] (2)

class #1:
real form #0: [0] (1)
real form #0: [1] (1)
\end{verbatim}

We compute
$$
\begin{aligned}
H^{-\tau}/(H^{-\tau})^0&=\{\diag(z,w,w,z)\,|\, z^2w^2=1\}/
\{\diag(ad,bc,bc,ad)\,|\, abcd=1\}\\
&=\{I, \diag(1,-1,-1,1)\}.
\end{aligned}
$$

Taking $z=I$ (this is {\tt class \#1}) we therefore have
\begin{equation}
\X_\delta(I)=\{\delta, \diag(1,-1,-1,1)\delta\}.
\end{equation}
These elements are not conjugate by $W$, so there are two strong real
forms mapping to  the
real form $SL(2,\mathbb H)$ ({\tt real form \#0}).

On the other hand take $z=-I$ ({\tt class \#0}). 
Then $x=\diag(1,1,-1,-1)\delta\in \X_\delta(-I)$, and 
$\X_{\delta}(-I)$ is obtained by multiplying x on the left
by $H^{-\tau}/(H^{-\tau})^0$. This gives
$$
\begin{aligned}
\X_{\delta}(-I)&=
\{\diag(1,1,-1,-1)\delta,\diag(1,-1,1,-1)\delta\}.
\end{aligned}
$$
These elements are conjugate by $W$, so there is only one strong real
form mapping to the real form $SL(4,\R)$
({\tt real form \#1}).
\end{exampleplain}

\bibliographystyle{amsalpha}

\enddocument